\documentclass{amsart}
\usepackage{amssymb}
\usepackage{latexsym}
\usepackage{amsthm}
\usepackage{amscd}
\theoremstyle{plain}

\newtheorem{thm}{Theorem}[section]
\newtheorem{cor}[thm]{Corollary}
\newtheorem{lem}[thm]{Lemma}
\newtheorem{prop}[thm]{Proposition}
\newtheorem{defn}{Definition}[section]

\numberwithin{equation}{section}

\newcommand{\im}{\mathrm{\ Im\ }}

\newcommand{\sgn}{\mathrm{sgn}}

\newcommand{\SL}{\mathit{SL}}
\newcommand{\iGamma}{\mathit\Gamma}
\newcommand{\tj}{\tilde{j}}
\newcommand{\tm}{\tilde{m}}
\newcommand{\tn}{{\tilde{n}}}
\newcommand{\new}{\mathrm{new}}
\newcommand{\old}{\mathrm{old}}

\newcommand{\lfl}{\left\lfloor}
\newcommand{\rfl}{\right\rfloor}

\newcommand{\BS}{\mathbf{S}}
\newcommand{\BT}{\mathbf{T}}
\newcommand{\ZZ}{\mathbb Z}
\newcommand{\CC}{\mathbb C}

\newcommand{\RR}{\mathbb R}

\newcommand{\ccsp}{\null\hskip 40pt}

\allowdisplaybreaks[1]

\begin{document}
\title[HECKE OPERATORS ON $\iGamma_0(2)$]
{Period polynomials and explicit formulas for Hecke operators on
  $\iGamma_0(2)$}
\author{Shinji Fukuhara}
\address{Department of Mathematics, Tsuda College, Tsuda-machi 2-1-1,
  Kodaira-shi, Tokyo 187-8577, Japan}
\email{fukuhara@tsuda.ac.jp}
\author{Yifan Yang}
\address{Department of Applied Mathematics, National Chiao Tung University,
  1001 Ta Hsueh Road, Hsinchu, Taiwan 300}
\email{yfyang@math.nctu.edu.tw}
\subjclass[2000]{Primary 11F25; Secondary 11F11, 11F67}
\keywords{Hecke operators, modular forms (one variable),
period polynomials}
\thanks{}

\begin{abstract} Let $S_{w+2}(\iGamma_0(N))$ be the vector space of
  cusp forms of weight $w+2$ on the congruence subgroup
  $\iGamma_0(N)$. We first determine explicit formulas for period
  polynomials of elements in $S_{w+2}(\iGamma_0(N))$ by means of
  Bernoulli polynomials. When $N=2$, from these explicit formulas we
  obtain new bases for $S_{w+2}(\iGamma_0(2))$, and extend the
  Eichler-Shimura-Manin isomorphism theorem to $\iGamma_0(2)$. This
  implies that there are natural correspondences between the spaces of
  cusp forms on $\iGamma_0(2)$ and the spaces of period polynomials.
  Based on these results, we will find explicit form of Hecke
  operators on $S_{w+2}(\iGamma_0(2))$. As an application of our main
  theorems, we will also give an affirmative answer to a speculation
  of Imamo\=glu and Kohnen on a basis of $S_{w+2}(\iGamma_0(2))$.
\end{abstract}

\maketitle

\begin{section}{Introduction}
Let $\iGamma$ be a congruence subgroup of $\SL_2(\ZZ)$. One of the
most important problems in the theory of modular forms is to
obtain explicit formulas for Hecke operators on cusp forms for
$\iGamma$. When $\iGamma$ is the full modular group $SL_2(\ZZ)$, this
was done in \cite{F5}, where we gave explicit formulas in terms of
Bernoulli numbers $B_k$ and divisor functions $\sigma_k$. Here we
briefly recall the approach in \cite{F5}.

For a cusp form $f$ of weight $w+2$ on $\SL_2(\ZZ)$ with $w\ge 2$ even
we consider the $n$-th period
$$
  r_n(f)=\int_0^{i\infty}f(z)z^n\,dz, \qquad 0\le n\le w.
$$
Since $r_n:S_{w+2}(\SL_2(\ZZ))\to\CC$ is a linear functional,
there exists a unique cusp form $R_n$ of weight $w+2$ such that
$r_n(f)=(f,R_n)$,
where $(f,g)$ is the Petersson inner product. In \cite{F5}, we first
showed that a certain subsets of $\{R_n\}$ forms a basis for the space
$S_{w+2}(\SL_2(\ZZ))$ of cusp forms of weight $w+2$. We then studied the
action of Hecke operators on this basis and obtained an explicit matrix
representation of Hecke operators. The Dedekind symbols (\cite{F4})
played a central role in our argument. (Note that periods and the Hecke
operators on periods have been studied by a number of mathematicians.
To name a few, see \cite{AN1,DI1,F5,KZ1,M1,M2,SH1,SK1,Z1,Z2}.)

The main goal in this article is to extend our formulas from
$\SL_2(\ZZ)$ to $\iGamma_0(N)$. Our starting point is the same as
the case for $\SL_2(\ZZ)$. Namely, for a cusp form $f$ of weight $w+2$
on $\iGamma_0(N)$, we consider the $n$-th period
$r_n(f)=\int^{i\infty}_0f(z)z^n\,dz$ for $0\le n\le w$, and let
$R_{\iGamma_0(N),w,n}$ be the unique cusp form determined by the
property
$$
  r_n(f)=2^{-1}(2i)^{w+1}(f,R_{\iGamma_0(N),w,n})
$$
for all cusp forms $f\in S_{w+2}(\iGamma_0(N))$.
(The choice of the constant $2^{-1}(2i)^{w+1}$ is to make our formulas
look nicer.) However, for lack of Dedekind symbols on $\iGamma_0(N)$,
the method for $\SL_2(\ZZ)$ in \cite{F5} would not work here, and we
ought to develop a new method to handle the case $\iGamma_0(N)$.

First of all, it is known in literature that these cusp forms
$R_{\iGamma_0(N),w,n}$ can be expressed in terms of certain Poincar\'e
series. We then observe that with those Poincar\'e series
representation, the action of Hecke operators on $R_{\iGamma_0(N),w,n}$
can be described very concretely. This observation turns out to be
crucial in our method. This is because if a basis in terms of
$R_{\iGamma_0(N),w,n}$ exists, then we can compute the periods of
$R_{\iGamma_0(N),w,n}$ and those of $T_m R_{\iGamma_0(N),w,n}$ to get
matrix representations of Hecke operators $T_m$. This basically
summarizes our approach. However, here we should remark that for a
general integer $N$, the set $\{R_{\iGamma_0(N),w,n}\}$ cannot
possibly span the whole space. This is because the dimension of
$S_{w+2}(\iGamma_0(N))$ is roughly $w[\SL_2(\ZZ):\iGamma_0(N)]/12$,
while we have only $w+1$ different $R_{\iGamma_0(N),w,n}$. For general
integers, one might need to consider more general periods. (See
\cite{AN1,DI1,SK1}.) Furthermore, for $m$ and $n$ with the same
parity, not much can be said about $r_m(R_{\iGamma_0(N),w,n})$. Thus,
for odd $n$, we consider only even periods, and for even $n$, we
consider only odd periods. Consequently, we will only consider the cases
where $[\SL_2(\ZZ):\iGamma_0(N)]\le 6$, that is, $N=2,\ldots,5$.

For the case $N=2$ we are able to determine a basis for every even
weight $w+2$. To achieve this, we show that the images of the
period maps are linearly independent for certain subsets of
$\{R_{\iGamma_0(2),w,n}\}$. Then we obtain explicit formulas for Hecke
operators by computing the periods of $T_mR_{\iGamma_0(2),w,n}$.
As a consequence of knowing an explicit basis for
$S_{w+2}(\iGamma_0(2))$, we can establish the $\iGamma_0(2)$-version
of the Eichler-Shimura-Manin theorem (\cite{E1}, \cite{KZ1},
\cite{M2}, \cite{SH1}), which strengthens a result of Imamo\=glu and
Kohnen (\cite[Proposition 3]{IK}. Furthermore, we give an affirmative
answer to a speculation raised by Imamo\=glu and Kohnen \cite{IK},
which is described as follows.

For even integers $k\ge 4$, let
$$
  E^{i\infty}_k(z)=\frac12\sum_{\gcd(c,d)=1,\,2|c}\frac1{(cz+d)^k}, \qquad
  E^0_k(z)=\frac12\sum_{\gcd(c,d)=1,\,2\nmid c}\frac1{(cz+d)^k}
$$
be normalized Eisenstein series of weight $k$ for the cusps $i\infty$
and $0$, respectively. It is obvious that
$E^{i\infty}_{2j}E^0_{k-2j}$, $j=2,\,3,\ldots,k/2-2$, are all
cusp forms of weight $k$ on $\iGamma_0(2)$. Inspired by a conjecture of
Chan and Chua \cite{CC}, Imamo\=glu and Kohnen \cite[Theorem 1]{IK}
proved that these forms generate the whole space $S_{k}(\iGamma_0(2))$
of cusp forms of weight $k$ on $\iGamma_0(2)$. The number of such
forms exceeds the dimension of $S_{k}(\iGamma_0(2))$ when $k>8$, so it
is natural to determine which functions form a basis for
$S_{k}(\iGamma_0(2))$. The experimental 
calculation in \cite{IK} suggests that the first few functions
$E_{2j}^{i\infty}E_{k-2j}^0$, $j=2,\,3,\ldots$, will
constitute a basis. In this paper, we will prove that this is indeed
true in general as a corollary of our main theorems.

For the cases $N=3,4,5$, we are unable to obtain analogous results,
although our numerical computation suggests that the first $d_w$
cusp forms $R_{\iGamma_0(N),w,2i}$, $i=1,\ldots,d_w=\dim
S_{w+2}(\iGamma_0(N))$, are always a basis for
$S_{w+2}(\iGamma_0(N))$. In any case, as long as a basis in terms of
$R_{\iGamma_0(N),w,n}$ is found, we can use our formulas for period
polynomials to obtain matrices for Hecke operators.
\end{section}

\begin{section}{Definitions and statements of results} Throughout the
  paper, we assume that $N$ is an integer with $N>1$ and $w$ is an
  even positive integer. For a cusp form $f\in S_{w+2}(\iGamma_0(N))$
  of weight $w+2$ on $\iGamma_0(N)$, we let
  $$
    r_n(f):=\int^{i\infty}_0f(z)z^n\,dz
  $$
  be the $n$-th period of $f$, and let the period polynomial $r(f)$ be
  defined by
  $$
    r(f)(X):=\int^{i\infty}_0f(z)(X-z)^w\,dz.
  $$
  Furthermore, even and odd period polynomials $r^+(f)$ and $r^-(f)$
  are defined by
  $$
    r^\pm(f)(X):=\frac12\{r(f)(X)\pm r(f)(-X)\}.
  $$
  If we set
  $$
    V_w:=\text{the vector space of polynomials of degree }\leq w
         \text{ in one variable }X,
  $$
  then we have homomorphisms of vector spaces:
  \begin{align*}
    r:\ &S_{w+2}(\iGamma_0(N))\to V_{w},  \\
    r^\pm:\ &S_{w+2}(\iGamma_0(N))\to V_{w}.
  \end{align*}

  Hereafter, $B_m(x)$ (resp. $B_m$) denotes the $m$-th Bernoulli
  polynomial (resp. number). By $B_m^0(x)$, we denote the $m$-th
  Bernoulli polynomial without its $B_1$-term (\cite[p.\ 208]{KZ1}):
  \begin{equation*}
    B_m^0(x)
      :=\sum_{\substack{0\leq i\leq m \\ i\ne 1}}
      \binom{m}{i}B_ix^{m-i}
     =\sum_{\substack{0\leq i\leq m \\ \text{$i$\ even}}}
      \binom{m}{i}B_ix^{m-i}.
  \end{equation*}
  Moreover, $\sgn(x)$ denotes the sign of $x\in\RR$. For an integer
  $n$ with $0\leq n\leq w$, let $\tn$ stand for $w-n$.

  First we need the following definitions to state our results:
  \begin{defn}\label{defn1.1}
  \begin{enumerate}
  \item For an integer $n$ such that $0\le n\le w$, we let
    $R_{\iGamma_0(N),w,n}(z)$ be the unique cusp form of weight
    $w+2$ on $\iGamma_0(N)$ characterized by 
    $$
      2^{-1}(2i)^{w+1}(f,R_{\iGamma_0(N),w,n})=\int_0^{i\infty}f(z)z^n\,dz
      =r_n(f)
    $$
    for all cusp forms $f$ of the same weight on $\iGamma_0(N)$, where
    $$
      (f,g)=\iint_{\iGamma_0(N)\backslash\mathbb H}
      f(z)\overline{g(z)}y^w\,dx\,dy, \qquad z=x+iy,
    $$
    denotes the Petersson inner product of $f$ and $g$;
  \item For an integer $n$ with $0<n<w$, we define a polynomial
    $S_{N,w,n}$ in $X$ by
    \begin{equation*}
      S_{N,w,n}(X):=
      \frac{N^\tn X^w}{\tn+1}B_{\tn+1}^0\left(\frac1{NX}\right)
      -\frac1{n+1}B_{n+1}^0(X).
    \end{equation*}
  \end{enumerate}
  \end{defn}
  The reader should be advised that our definition of the Petersson
  inner product differs from the standard definition by a factor
  $[SL_2(\ZZ):\iGamma_0(N)]$. This is for the sake of making the
  presentations simpler.

  In our first theorem we evaluate the period polynomials of
  $R_{\iGamma_0(N),w,n}$. This result is crucial in our subsequent
  discussion.

  \begin{thm}\label{thm1.1} Let $N$ be an integer greater than $1$.
  For an even integer $n$ with $0<n<w$, we have
  \begin{equation*}
    r^-(R_{\iGamma_0(N),w,n})(X)=S_{N,w,n}(X).
  \end{equation*}
  Also, for an odd integer $n$ with $0<n<w$, we have
  \begin{equation*}
  \begin{split}
   &r^+(R_{\iGamma_0(N),w,n})(X)=S_{N,w,n}(X) \\
   &\qquad-\frac{(w+2)B_{n+1}B_{\tn+1}}{(n+1)(\tn+1)B_{w+2}}
    \left(\frac{X^w}N\prod_{p|N}\frac{1-p^{-(n+1)}}{1-p^{-(w+2)}}
   -\frac1{N^{n+1}}\prod_{p|N}\frac{1-p^{-(\tn+1)}}{1-p^{-(w+2)}}\right),
  \end{split}
  \end{equation*}
  where $p$ runs over all prime divisors of $N$.
  \end{thm}

  We remark that the formulas for period polynomials for the case
  $N>1$ is actually simpler than those for the case $N=1$ given in
  \cite{KZ1}. This is because matrices in $\iGamma_0(N)$ with $N>1$
  cannot have zeroes at the $(1,1)$-entry and the $(2,2)$-entry. We
  note also that the formula for the period polynomial
  $r^+(R_{\iGamma_0(N),w,n})(X)$, $n$ odd, has already been discussed
  in \cite{AN1}. However, the discussion in \cite{AN1} is in a
  very broad setting, and consequently, the formulas of \cite{AN1} are
  too complicated to be applied readily to our situation.

  We now devote our attention to the study of the action of Hecke
  operators on $R_{\iGamma_0(N),w,n}$. For this purpose, we will need
  the following definitions:

  \begin{defn}\label{defn1.2}
  \begin{enumerate}
  \item For a positive integer $m$, let
    \begin{equation*}
      H_{N,m}:=\left\{
        \begin{bmatrix}a&b\\c&d\end{bmatrix}
        \Big|\ a,b,c,d\in\ZZ;\ ad-bc=m;\
            c\equiv 0 \pmod {N};\ \gcd(a,N)=1
        \right\};
    \end{equation*}
  \item For positive integers $m$ and $n$ such that $0<n<w$, let
  \begin{align*}
    &R_{\iGamma_0(N),w,n}^m(z) \\
    &\qquad:=m^{w+1}c_{w,n}^{-1}
      \sum_{\substack{\left[\begin{smallmatrix}a&b\\c&d\end{smallmatrix}\right]
                      \in H_{N,m}}}
           \frac{1}{(az+b)^{\tn+1}(cz+d)^{n+1}},
      \ c_{w,n}=(-1)^n2\pi i\binom wn;
  \end{align*}
  \item For positive integers $m$ and $n$ such that $0<n<w$, we define
    a polynomial $S_{N,w,n}^m$ in $X$ by
    {\allowdisplaybreaks
    \begin{align*}
      S_{N,w,n}^m(&X) :=\frac{1}{2}
      \sum_{\substack{\left[\begin{smallmatrix}a&b\\c&d\end{smallmatrix}\right]
                      \in H_{N,m} \\
                      abcb<0}}
      \sgn(ab)(aX+b)^{n}(cX+d)^{\tn} \\
      &\ccsp+\sum_{\substack{ad=m,\ a>0 \\ \gcd(a,N)=1}}
        \left\{\frac{a^nN^{\tn}X^w}{\tn+1}B^0_{\tn+1}(\frac{d}{NX})
        -\frac{d^{\tn}}{n+1}B^0_{n+1}(aX)\right\}.
    \end{align*}
    }
  \end{enumerate}
  \end{defn}

  Then the actions of Hecke operators on
  $R_{\iGamma_0(N),w,n}$ are explicitly as follows:
  \begin{lem}\label{thm1.5} Let $m$ be a positive integer, the action
  of the Hecke operators $T_m$ on $R_{\iGamma_0(N),w,n}$ is
  \begin{equation*}
    T_m(R_{\iGamma_0(N),w,n})=R_{\iGamma_0(N),w,n}^m.
  \end{equation*}
  \end{lem}

  We remark that Lemma \ref{thm1.5} for the case $N=1$ is casually
  mentioned in \cite{KZ1} without a proof. For the convenience of the
  reader we will provide a proof in Section \ref{sect3}.

  From this series expression of $T_m(R_{\iGamma_0(N),w,n})$, we can
  compute the period polynomials of $T_m(R_{\iGamma_0(N),w,n})$ for
  even $n$.
  \begin{thm}\label{thm1.6}
  Let $n$ be an even integer with $0<n<w$. If $m$ is a positive
  integer not divisible by $N$, then
  \begin{equation*}
    r^-(R_{\iGamma_0(N),w,n}^m)(X)=S_{N,w,n}^m(X).
  \end{equation*}
  Furthermore, when $m$ is a multiple of $N$, we have
  \begin{equation*}
    r^-(R_{\iGamma_0(N),w,n}^m)(X)=S_{N,w,n}^m(X)
   -\frac{(NX)^w}{n+1}\sum_{d|N}\frac{\mu(N/d)}{d^n}\sum_{c|(m/N)}c^\tn
    B_{n+1}^0\left(\frac{md}{cN^2X}\right),
  \end{equation*}
  where $\mu(d)$ is the M\"obius function.
  \end{thm}

  \noindent{\bf Examples.}
  \begin{enumerate}
  \item Let $N=4$, $w=4$ and consider the Hecke operator $T_2$ on
  $S_6(\iGamma_0(4))$. Since $T_2$ sends a cusp form on $\iGamma_0(4)$
  to a cusp form on $\iGamma_0(2)$ and $\dim
  S_6(\iGamma_0(2))=0$, we should have $S_{4,4,2}^2(X)=0$. Indeed, we
  have
  \begin{equation*}
  \begin{split}
    S_{4,4,2}^2(X)&=\frac{4^2X^4}3B_3^0\left(\frac2{4X}\right)
   -\frac{2^2}3B_3^0(X) \\
    &=\frac{16X^4}3\left(\frac1{8X^3}+\frac36
    \frac1{2X}\right)
   -\frac43\left(X^3+\frac36X\right)=0.
  \end{split}
  \end{equation*}
  \item Let $N=4$, $w=6$, and consider the Hecke operator
  $T_8=T_2^3$. Then $T_8R_{\iGamma_0(4),6,2}$ is a scalar multiple of
  $R_{\iGamma_0(2),6,2}$ since $\dim S_8(\iGamma_0(2))=1$. Thus, the
  period polynomial $r^-(R^8_{\iGamma_0(4),6,2})(X)$ should be a
  scalar multiple of $r^-(R_{\iGamma_0(2),6,2})(X)$. We have
  \begin{equation*}
  \begin{split}
    S_{2,6,2}(X)=\frac{2^4X^6}5B_5^0\left(\frac1{2X}\right)
   -\frac13B_3^0(X)=-\frac1{15}(4X^5-5X^3+X).
  \end{split}
  \end{equation*}
  For $r^-(R^8_{\iGamma_0(4),6,2})$, we find that there are $4$ tuples
  $(a,b,c,d)=\pm(1,1,-4,4)$ and $\pm(1,-1,4,4)$ contributing to the first
  sum in $S_{4,6,2}^8(X)$, giving
  $$
    -1024(X^5-2X^3+X)
  $$
  while the remaining part of $S_{4,6,2}^8(X)$ is
  $$
    \frac{4^4X^6}5B_5^0(2/X)-\frac{8^4}3B_3^0(X)
   =-\frac{256}{15}(X^5+40X^3-56X).
  $$
  The other term in $r^-(R^8_{\iGamma_0(4),6,2})$ is
  \begin{equation*}
  \begin{split}
   -\frac{(4X)^6}3\sum_{d|4}\mu(4/d)d^{-2}\sum_{c|2}c^4B_3^0(d/(2cX))
   =256(3X^5-4X^3),
  \end{split}
  \end{equation*}
  and we find
  $$
    r^-(R^8_{\iGamma_0(4),6,2})(X)
   =-\frac{1024}{15}(4X^5-5X^3+X),
  $$
  which is indeed a multiple of $r^-(R_{\iGamma_0(2),6,2})(X)$.
  \end{enumerate}

  \medskip

  \noindent{\bf Remark.} Observant readers will notice that when the
  space $S_{w+2}(\iGamma_0(N))$ has dimension $1$, we have 
  $$
    R^m_{\iGamma_0(N),w,n}=a_m R_{\iGamma_0(N),w,n},
  $$
  where $a_m$ is the eigenvalue of the $m$-th Hecke operator.
  Thus, by considering the coefficients of the period polynomials of
  both sides, we can obtain explicit expressions for $a_m$. For
  example, take $N=2$, $w=6$. If we choose $n=2$, then Theorem
  \ref{thm1.1} yields
  $$
    r^-(R_{\iGamma_0(2),6,2})(X)=N^3B_4X^5+\cdots,
  $$
  while Lemma \ref{thm1.5} and Theorem \ref{thm1.6} give, for $2\nmid
  m$,
  \begin{equation*}
  \begin{split}
    r^-(T_mR_{\iGamma_0(2),6,2})&=X^5N^3B_4\sum_{ad=m}a^2d \\
   &\qquad+2X^5\sum_{\substack{a,b,c,d>0\\ ad+2bc=m}}\left(
    2ab(2c)^4-4a^2(2c)^3d\right)+\cdots.
  \end{split}
  \end{equation*}
  Thus, for $2\nmid m$, we see that the $m$-th Hecke eigenvalue, is
  $$
    m\sum_{a|m}a+240\sum_{\substack{a,b,c,d>0\\ ad+2bc=m}}
    ac^3(ad-bc)
   =m\sigma_1(m)+240\sum_{u+2v=m}(u-v)\sigma_1(u)\sigma_3(v).
  $$
  Although this formula looks like an interesting application of
  period polynomials, it, in fact, can be more easily obtained by
  applying the Rankin-Cohen bracket (see \cite{Cohen,Rankin}) to the
  pair of modular forms
  $$
    1-24\sum_{n=1}^\infty n\left(\frac{2q^{2n}}{1-q^{2n}}
     -\frac{q^n}{1-q^n}\right), \qquad
    1+240\sum_{n=1}^\infty\frac{n^3q^{2n}}{1-q^{2n}}, \qquad
    (q=e^{2\pi iz})
  $$
  on $\iGamma_0(2)$.
  \medskip

  In the remainder of the section, we restrict our attention to the
  case $N=2$. In this case we are able to determine bases for
  $S_{w+2}(\iGamma_0(2))$. For the sake of convenience, we set
  \begin{equation*}
    d_w:=\lfl\frac{w-2}{4}\rfl\
    \text{(i.e., $d_w$ is the greatest integer not exceeding
    $(w-2)/4$).}
  \end{equation*}
  It is well-known (see e.g. \cite{DS1}) that
  \begin{equation*}
        \dim S_{w+2}(\iGamma_0(2))=d_w.
  \end{equation*}
  By applying Theorem \ref{thm1.1} with $N=2$ and showing the linear
  independence among the period polynomials of $R_{\iGamma_0(2),w,n}$
  for certain ranges of integers $n$, we can determine a basis for
  $S_{w+2}(\iGamma_0(2))$.

  \begin{thm}\label{thm1.2} Each of the sets
  \begin{equation*}
    \{R_{\iGamma_0(2),w,{2i}}\ |\ i=1,2,\ldots,d_w \}, \qquad
    \{R_{\iGamma_0(2),w,{w-2i}}\ |\ i=1,2,\ldots,d_w\}
  \end{equation*}
  and
  \begin{equation*}
    \{R_{\iGamma_0(2),w,{2i-1}}\ |\ i=1,2,\ldots,d_w \}, \qquad
    \{R_{\iGamma_0(2),w,{w-2i+1}}\ |\ i=1,2,\ldots,d_w\}
  \end{equation*}
  forms a basis for $S_{w+2}(\iGamma_0(2))$.
  \end{thm}

  Now if $f$ is a cusp form of weight $w+2$ on $\iGamma_0(2)$ such
  that $r_{2i-1}(f)=0$ for all $i=1,\ldots,d_w$, then we have
  $(f,R_{\iGamma_0(2),w,{2i-1}})=0$ for all $i=1,\ldots,d_w$. In view
  of Theorem \ref{thm1.2}, this implies that $f$ is identically zero.
  Thus, we see that the map $f\mapsto
  (r_1(f),r_3(f),\ldots,r_{2d_w-1}(f))$ is an isomorphism and
  consequently $r^-:S_{w+2}(\iGamma_0(2))\to V_w$ is injective. The
  same argument applies to other similar maps as well. This is the
  $\iGamma_0(2)$-version of the Eichler-Shimura-Manin theorem
  (\cite{E1}, \cite{KZ1}, \cite{M2}, \cite{SH1}).

  \begin{cor}\label{cor1.3}
  The maps
  \begin{equation*}
    r^+,\,r^-:\ S_{w+2}(\iGamma_0(2))\to V_{w}
  \end{equation*}
  are both injective.
  \end{cor}

  Note that the $r^-$ case in Corollary \ref{cor1.3} was proved in
  \cite{IK} by a different method. However, our result is stronger in
  the sense that we have actually shown that the function
  $$
    f\mapsto(r_1(f),\ldots,r_{2d_w-1}(f))
  $$
  and other similar functions are isomorphisms between the vector
  spaces $S_{w+2}(\iGamma_0(2))$ and $\CC^{d_w}$.

  The proof of Theorem \ref{thm1.2} can be modified to yield
  alternative bases for $\iGamma_0(2)$, which in turn implies the
  speculation of Imamo\=glu and Kohnen \cite{IK} mentioned in the
  introduction that the first $d_w$ products
  $E_{2j+2}^{i\infty}E_{w-2j}^0$, $j=1,\ldots,d_w$, form a basis.

  \begin{thm} \label{thm1.4}
    Let $E^{i\infty}_{2j}$ and $E^0_{2j}$ be the normalized
    Eisenstein series of weight $2j$ on $\iGamma_0(2)$ associated with
    the cusps $i\infty$ and $0$, respectively. Then each of the two
    sets
    $$
      \{E^0_{2j+2}E^{i\infty}_{w-2j}\ |\  j=1,\ldots,d_w\}, \quad
      \{E^0_{w-2j}E^{i\infty}_{2j+2}\ |\  j=1,\ldots,d_w\}
    $$
    forms a basis for $S_{w+2}(\iGamma_0(2))$.
  \end{thm}
  The proof of Theorem \ref{thm1.4} is slightly more involved than
  that of Theorem \ref{thm1.2} since it requires an exact evaluation of
  Hankel determinants of Bernoulli numbers. The proof of the
  determinant formulas is supplied by Christian Krattenthaler, and
  will be given in Appendix A.

  Finally, as an application of Theorems \ref{thm1.1},
  \ref{thm1.6} and \ref{thm1.2}, we will give explicit formulas for
  the Hecke operators on $S_{w+2}(\iGamma_0(2))$. Let
  \begin{equation*}
    f(X)=\sum_{\nu=0}^{w}a_\nu X^{w-\nu}
    \text{\ \ \ and\ \ \ }
    g(X)=\sum_{\nu=0}^{w}b_\nu X^{w-\nu}
  \end{equation*}
  be polynomials in $X$ with degree $\leq w$. Then their inner product
  $\langle f,g\rangle$ is defined by
  \begin{equation*}
    \langle f,g\rangle:=\sum_{\nu=0}^{w}a_\nu \bar{b}_\nu,
  \end{equation*}
  where $\bar{b}_\nu$ denotes the complex conjugate of $b_\nu$.

  Under this notation we obtain the following result.
  \begin{thm}\label{thm1.7} Let $m$ be a positive integer, and let
    $\BT_{m}$ be the matrix representing the Hecke operator
    \begin{equation*}
      T_m:S_{w+2}(\iGamma_0(2))\to S_{w+2}(\iGamma_0(2))
    \end{equation*}
    with respect to the basis
    \begin{equation*}
      c_{w,2i}R_{\iGamma_0(2),w,{2i}}\ \ (i=1,2,\ldots,d_w).
    \end{equation*}
  Let $\BS_1$ and $\BS_2$ be matrices defined by
  \begin{equation*}
    \BS_1:=
    \begin{bmatrix}
    \langle S_{2,w,{2i}},S_{2,w,{2j}}\rangle
    \end{bmatrix}
    \ \ \ (i,j=1,2,\ldots,d_w),
  \end{equation*}
  \begin{equation*}
    \BS_2:=
     \begin{bmatrix}
     \langle S_{2,w,{2i}},S_{2,w,{2j}}^m\rangle
     \end{bmatrix}
     \ \ \ (i,j=1,2,\ldots,d_w).
  \end{equation*}
  Then $\BT_{m}$ can be expressed as
  \begin{equation*}
    \BT_{m}=\BS_1^{-1}\BS_2.
  \end{equation*}
  \end{thm}
  Consequently, this gives an algorithm for computing the matrix
  $\BT_{m}$ representing the Hecke operator $T_{m}$ on
  $S_{w+2}(\iGamma_0(2))$. We append a computer program in the last
  section.

  We conclude this section with two examples.
  \medskip

  \noindent{\bf Examples.}
  \begin{enumerate}
  \item Consider the Hecke operator $T_2$ on
  $S_{12}(\iGamma_0(2))$, which is of dimension $2$. We have
  \begin{align*}
    S_{2,10,2}(X)&=-\frac1{45}(192X^9-320X^7+168X^5-45X^3+5X), \\
    S_{2,10,4}(X)&=\frac1{210}(160X^9-280X^6+168X^4-55X^3+7X), \\
    S_{2,10,2}^2(X)&=\frac{128}{45}(12X^9+5X^7-42X^5+30X^3-5X), \\
    S_{2,10,4}^2(X)&=-\frac{32}{105}(44X^9-35X^7-42X^5+40X^3-7X).
  \end{align*}
  Thus,
  $$
    \BS_1^{-1}\BS_2=\begin{pmatrix}
      -208 & 36 \\ -1120 & 184 \end{pmatrix},
  $$
  whose characteristic polynomial is $x^2+24x+2048$. To check the
  correctness, we observe that the space $S_{12}(\iGamma_0(2))$ is
  spanned by $\Delta(z)$ and $\Delta(2z)$, where
  $$
    \Delta(z)=e^{2\pi iz}\prod_{n=1}^\infty(1-e^{2\pi inz})^{24}.
  $$
  We find
  $$
    T_2\Delta(z)=-24\Delta(z)-2048\Delta(2z), \qquad
    T_2\Delta(2z)=\Delta(z).
  $$
  It is clear that the characteristic polynomial of $T_2$ is indeed
  $x^2+24x+2048$.
  \item Consider $N=4$ and $w=8$. The space $S_{10}(\iGamma_0(4))$ has
  dimension 3 with a two-dimensional oldspace coming from
  $S_{10}(\iGamma_0(2))$ and a one-dimensional newspace. We can show
  that $R_{\iGamma_0(4),8,2i}$, $i=1,2,3$, form a basis. Thus, by the
  same reasoning as that for Theorem \ref{thm1.7}, the matrix of a
  Hecke operator $T_m$ with respect to $c_{w,2i}R_{\iGamma_0(4),8,2i}$
  would be $\BS_1^{-1}\BS_2$, where
  \begin{align*}
    \BS_1&=[\langle S_{4,8,2i},S_{4,8,2j}\rangle] \quad
    (i,\,j=1,2,3), \\
    \BS_2&=[\langle S_{4,8,2i},S_{4,8,2j}^m\rangle] \quad
    (i,\,j=1,2,3).
  \end{align*}
  Take $m=3$, for example. We find
  $$
    \BS_1^{-1}\BS_2=\frac1{152915}\begin{pmatrix}
    2456678965260 & -224610211392 & 61847064000 \\
    37961609400000 & -3470759119380 & 955676880000 \\
    40281954570000 & -3682878636192 & 1014067309260
    \end{pmatrix},
  $$
  whose characteristic polynomial is $(x-228)(x+156)^2$. Therefore,
  the eigenvalues of $T_3$ of the Hecke eigenform on $\iGamma_0(2)$
  and the newform on $\iGamma_0(4)$ of weight $10$ are $-156$ and
  $228$, respectively. Indeed, the Hecke eigenform on $\iGamma_0(2)$
  of weight $10$ is
  $$
    \eta(z)^8\eta(2z)^8\left(1-24\sum_{n=1}^\infty
    \left(\frac{2nq^{2n}}{1-q^{2n}}-\frac{nq^n}{1-q^n}\right)\right)
   =q+16q^2-156q^3+256q^4+\cdots,
  $$
  and the newform on $\iGamma_0(4)$ of weight $10$ is
  $$
    \eta(2z)^{12}\left(1+240\sum_{n=1}^\infty\frac{n^3q^{2n}}{1-q^{2n}}
    \right)=q+228q^3-666q^5+\cdots.
  $$
  \end{enumerate}
\end{section}

\begin{section}{Proof of Theorem \ref{thm1.1}} \label{sect2}
  First we give a characterization of the cusp form
  $R_{\iGamma_0(N),w,n}(z)$ in terms of certain Poincar\'{e} series.
  Recall that $(f,g)$ denotes the Petersson inner product
  of $f$ and $g$ in $S_{w+2}(\iGamma_0(N))$. For
  $\gamma=\begin{bmatrix}a&b\\c&d\end{bmatrix}$ in $GL^+_2(\RR)$,
  \ $f|_\gamma$ is defined by
  \begin{equation*}
  (f|_\gamma)(z):=(\det\gamma)^{w/2+1}(cz+d)^{-w-2}f(\frac{az+b}{cz+d}).
  \end{equation*}

  \begin{lem}\label{lem2.1}
    For $0<n<w$, the cusp form $R_{\iGamma_0(N),w,n}$ is given by
    \begin{equation} \label{eqn: lem2.1}
      R_{\iGamma_0(N),w,n}=c_{w,n}^{-1}\sum_{\substack{\left[
        \begin{smallmatrix}a&b\\ c&d\end{smallmatrix}\right]
        \in\iGamma_0(N)}}
        \frac{1}{(az+b)^{\tn+1}(cz+d)^{n+1}}, \quad
      c_{w,n}=(-1)^n2\pi i\binom wn.
    \end{equation}
  \end{lem}

  \begin{proof} The proof for the case $N=1$ is given in
    \cite{C1}. (See also \cite{KZ1}.)  For general $N$, the lemma is
    just a special case of Proposition 3 of \cite{AN1}. We reproduce
    it here for the sake of completeness.

    Let $R_n$ denote the function on the right-hand side of
    \eqref{eqn: lem2.1}. We have
    \begin{equation*}
    \begin{split}
    (f,c_{w,n}R_n)
  &=\left(f,\sum_{\gamma\in\iGamma_0(N)}z^{-\tn-1}|_\gamma\right)
   =2\int_0^\infty y^w\left(\int^\infty_{-\infty}
    \frac{f(x+iy)}{(x-iy)^{\tn+1}}\,dx\right)\,dy \\
  &=2\int^\infty_0 y^w\left(\frac{2\pi i}{\tn!}f^{(\tn)}(2iy)\right)\,dy \\
  &=4\pi i\binom wn\left(\frac i2\right)^\tn\int^\infty_0
    y^{w-\tn}f(2iy)\,dy \\
  &=4\pi i\binom wn\left(\frac i2\right)^\tn(2i)^{-n-1}r_n(f).
    \end{split}
    \end{equation*}
    Then the lemma follows.
  \end{proof}

  Next we study how $R_{\iGamma_0(N),w,n}(z)$ behaves under the
  Atkin-Lehner involution $\left[\begin{smallmatrix}0&-1\\
  N & 0\end{smallmatrix}\right]$:

  \begin{lem}\label{lem2.2} For all integers $m$ and $n$ with $0\le
  m\le w$ and $0<n<w$, we have
  $$
    R_{\iGamma_0(N),w,n}(\frac{-1}{Nz})=(-N)^{\tn+1}z^{w+2}
    R_{\iGamma_0(N),w,\tn}(z)
  $$
  and
  $$
    r_m(R_{\iGamma_0(N),w,n})=(-N)^{\tn-m}r_{\tm}(R_{\iGamma_0(N),w,\tn}).
  $$
  \end{lem}

  \begin{proof}
    Recall that the Atkin-Lehner involution $\omega_N$ defined by
    $$
      f\big|_{\omega_N}(z)=(\sqrt N z)^{-w-2}f(\frac{-1}{Nz})
    $$
    is self-adjoint with respect to the Petersson inner product. Thus,
    we have, for all cusp forms $f$ of weight $w+2$ on $\iGamma_0(N)$,
    $$
      (f,R_{\iGamma_0(N),w,n}|_{\omega_N})
     =(f|_{\omega_N},R_{\iGamma_0(N),w,n})
     =2(2i)^{-w-1}\int_0^{i\infty}(\sqrt Nz)^{-w-2}f(\frac{-1}{Nz})z^n\,dz.
    $$
    We then make a change of variable $z\mapsto-1/Nz$ in the integral.
    We find that
    \begin{equation*}
    \begin{split}
      (f,R_{\iGamma_0(N),w,n}|_{\omega_N})
    &=2(2i)^{-w-1}\int^0_{i\infty}(\sqrt Nz)^{w+2}
      f(z)(-Nz)^{-n}\frac{dz}{Nz^2} \\ 
    &=(-1)^{n+1}N^{w/2-n}(f,R_{\iGamma_0(N),w,\tn}).
    \end{split}
    \end{equation*}
    It follows that
    $$
      R_{\iGamma_0(N),w,n}\big|_{\omega_N}=(-1)^{n+1}N^{w/2-n}
      R_{\iGamma_0(N),w,\tn},
    $$
    which is equivalent to our first assertion. We now prove the
    second assertion.

    Since $w_N$ is an isometry with respect to the Petersson inner
    product, we have
    \begin{equation*}
    \begin{split}
     2(2i)^{-w-1}r_m(R_{\iGamma_0(N),w,n})
    &=(R_{\iGamma_0(N),w,n},R_{\iGamma_0(N),w,m}) \\
    &=(R_{\iGamma_0(N),w,n}|_{\omega_N},R_{\iGamma_0(N),w,m}|_{\omega_N}).
    \end{split}
    \end{equation*}
    Then by the first part of the lemma, we see that
    \begin{equation*}
    \begin{split}
      r_m(R_{\iGamma_0(N),w,n})&=2^{-1}(2i)^{w+1}(-1)^{m+n}N^{w-n-m}
      (R_{\iGamma_0(N),w,\tn},R_{\iGamma_0(N),w,\tm}) \\
      &=(-N)^{\tn-m}r_{\tm}(R_{\iGamma_0(N),w,\tn}).
    \end{split}
    \end{equation*}
    This completes the proof of the lemma.
  \end{proof}

  Now we are ready to prove the following proposition which is a key
  for our arguments:

  \begin{prop}\label{prop2.3} Let $c_{w,n}=(-1)^n2\pi i\binom wn$.
    For $0\le m\le w$ and $0<n<w$ with opposite parity, we have
  \begin{enumerate}
  \item if $m+n>w$,
  \begin{equation*}
  \begin{split}
    c_{w,n}r_m(R_{\iGamma_0(N),w,n})&=\frac{2\pi i}{m+1}\binom{m+1}\tn
    B_{m-\tn+1}-\frac{2\pi i}{Nn}\delta_{m,\tn+1} \\
   &\qquad-\delta_{m,w}\binom{w+2}{n+1}
    \frac{2\pi iB_{n+1}B_{\tn+1}}{(w+1)N^{n+1}B_{w+2}}
    \prod_{p|N}\frac{1-p^{-(\tn+1)}}{1-p^{-(w+2)}};
  \end{split}
  \end{equation*}
  \item if $m+n<w$,
  \begin{equation*}
  \begin{split}
    c_{w,n}r_m(R_{\iGamma_0(N),w,n})&=(-N)^{\tn-m}
      \Bigg(\frac{2\pi i}{\tm+1}\binom{\tm+1}nB_{\tm-n+1}
     -\frac{2\pi i}{N\tn}\delta_{\tm,n+1} \\
     &\qquad-\delta_{m,0}\binom{w+2}{n+1}
      \frac{2\pi iB_{n+1}B_{\tn+1}}{(w+1)N^{\tn+1}B_{w+2}}
      \prod_{p|N}\frac{1-p^{-(n+1)}}{1-p^{-(w+2)}}\Bigg),
  \end{split}
  \end{equation*}
  where $\delta_{i,j}$ is the Kronecker delta symbol.
  \end{enumerate}

  Furthermore, for an odd integer $m$ with $0<m<w$, we have
  \begin{equation*}
  \begin{split}
    r_m(R_{\iGamma_0(N),w,w})
  &=\frac{B_{m+1}}{m+1}-\frac{\delta_{w,\tm+1}}{Nw} \\
   &\qquad -\frac{(w+2)B_{m+1}B_{\tm+1}}{N^{m+1}(m+1)(\tm+1)B_{w+2}}
    \prod_{p|N}\frac{1-p^{-(\tm+1)}}{1-p^{-(w+2)}}
  \end{split}
  \end{equation*}
  and
  \begin{equation*}
  \begin{split}
    r_m(R_{\iGamma_0(N),w,0})
  &=-N^{\tm}\Bigg(\frac{B_{\tm+1}}{\tm+1}-\frac{\delta_{w,m+1}}{Nw} \\
  &\qquad -\frac{(w+2)B_{m+1}B_{\tm+1}}{N^{\tm+1}(m+1)(\tm+1)B_{w+2}}
   \prod_{p|N}\frac{1-p^{-(m+1)}}{1-p^{-(w+2)}}\Bigg).
  \end{split}
  \end{equation*}
  \end{prop}

  \begin{proof} To ease the notations, throughout the proof, we
  write
  $$
    R_n=R_{\iGamma_0(N),w,n}=c^{-1}_{w,n}
    \sum_{\substack{\left[\begin{smallmatrix}a&b\\c&d\end{smallmatrix}\right]
                      \in \iGamma_0(N)}}
    \frac{1}{(az+b)^{\tn+1}(cz+d)^{n+1}}.
  $$
  Here we first consider the cases where $m>\tn>0$ (that is, $m+n>w$);
  the results for the remaining cases will follow from these
  cases using Lemma \ref{lem2.2}. Now for such $m$ and $n$, we partition
  the sum defining $R_n$ into three parts,
  $$
    c_{w,n}R_n(z)=\sum_{c=0}+\sum_{b=0,c\neq 0}+\sum_{bc\neq 0},
  $$
  and let $S_i$, $i=1,2,3$, denote the $i$-th sum, respectively.

  For $S_1$ we have
  \begin{equation} \label{S1}
    S_1(z)=\sum_{\substack{
    \left[\begin{smallmatrix}1&b\\ 0&1\end{smallmatrix}\right]}}
     +\sum_{\substack{
    \left[\begin{smallmatrix}-1&-b\\ 0&-1\end{smallmatrix}\right]}}
     =2\sum_{b\in\ZZ}\frac1{(z+b)^{\tn+1}}
   =\frac{2(-2\pi i)^{\tn+1}}{\tn !}
    \sum_{r=1}^\infty r^{\tn}e^{2\pi irz},
  \end{equation}
  where in the last equality we use the formula
  $$
    \sum_{b\in\ZZ}(z+b)^{-n}=\frac{(-2\pi i)^n}{\Gamma(n)}
    \sum_{r=1}^\infty r^{n-1}e^{2\pi irz},
  $$
  valid for all positive integers $n$ and complex numbers $z$ with
  $\im z>0$. It follows that
  \begin{equation*}
  \begin{split}
    \int_0^{i\infty}z^m S_1(z)\,dz
  &=\frac{2(-2\pi i)^{\tn+1}}{\tn!}\int^{i\infty}_0 z^m
    \sum_{r=1}^\infty r^\tn e^{2\pi irz}\,dz \\
  &=2(-2\pi i)^{\tn-m}\frac{m!}{\tn!}\sum_{r=1}^\infty r^{-(m-\tn+1)}
   =2(-2\pi i)^{\tn-m}\frac{m!}{\tn!}\zeta(m-\tn+1).
  \end{split}
  \end{equation*}
  Using the well-known formula
  \begin{equation}\label{eqn2.1}
    \zeta(2p)=\frac{(-1)^{p-1}}{(2p)!}2^{2p-1}\pi^{2p}B_{2p} \ \ (p>0),
  \end{equation}
  we simplify the last expression to
  \begin{equation}\label{eqn2.2}
    \int_0^{i\infty}z^m S_1(z)\,dz
   =\frac{2\pi i}{m+1}\binom{m+1}\tn B_{m-\tn+1}.
  \end{equation}

  We next consider the contribution from
  $$
    S_2(z)=2\sum_{c\in\ZZ,c\neq 0}\frac1{z^{\tn+1}(cNz+1)^{n+1}}.
  $$
  We first show that the integral
  $$
    \int^{i\infty}_0z^m S_2(z)\,dz
  $$
  converges absolutely when $m>\tn+1$. For $z$ with $|z|\ll 1$ we have
  $$
    \frac 1{|cNz+1|^{n+1}}\ll
    \begin{cases}
      1, &\text{if }|c|\ll 1/|z|, \\
      |cz|^{-n-1}, &\text{if }|c|\gg 1/|z|,
    \end{cases}
  $$
  and hence
  $$
    \sum_{c\in\ZZ,c\neq 0}\frac1{|cNz+1|^{n+1}}
    \ll\sum_{|c|\ll 1/|z|}1+\frac1{|z|^{n+1}}\sum_{|c|\gg 1/|z|}
    \frac 1{|c|^{n+1}}\ll\frac 1{|z|}.
  $$
  For $|z|\gg 1$ we have
  $$
    \frac 1{|cNz+1|^{n+1}}\ll\frac1{|cz|^{n+1}}
  $$
  and
  \begin{equation}\label{eqn2.3}
    \sum_{c\in\ZZ,c\neq 0}\frac1{|cNz+1|^{n+1}}\ll\frac 1{|z|^{n+1}}.
  \end{equation}
  It follows that, for $m>\tn+1$,
  $$
    \int_0^{i\infty}|z|^m\sum_{c\in\ZZ,c\neq 0}
      \frac1{|z|^{\tn+1}|cNz+1|^{n+1}}\,dz
    \ll\int_0^{1}t^{m-\tn-2}\,dt+
      \int_1^\infty t^{m-\tn-n-2}\,dt<\infty.
  $$
  Now, having proved that $\int_0^{i\infty}z^m S_2(z)\,dz$
  converges absolutely for $m>\tn+1$, we can change the order of
  integration and summation freely at will. We have
  \begin{equation*}
  \begin{split}
    \int_0^{i\infty}z^mS_2(z)\,dz
  &=2\int_0^{i\infty}z^{m-\tn-1}\sum_{c\in\ZZ,c\neq 0}
    \frac1{(cNz+1)^{n+1}}\,dz \\
  &=\sum_{c\in\ZZ,c\neq 0}\left\{\int_0^{i\infty}
    \frac{z^{m-\tn-1}}{(cNz+1)^{n+1}}\,dz
   +\int^0_{-i\infty}
    \frac{(-z)^{m-\tn-1}}{(-cNz+1)^{n+1}}\,dz\right\} \\
  &=\sum_{c\in\ZZ,c\neq 0}\int^{i\infty}_{-i\infty}z^{m-\tn-1}
    \frac1{(cNz+1)^{n+1}}\,dz.
  \end{split}
  \end{equation*}
  We then move the line of integration to the far left or the far right
  according to the sign of $c$. Thus, the integral is $0$, provided that
  $m>\tn+1$.

  For the case $m=\tn+1$, we first observe that, for large $z$, we have
  \eqref{eqn2.3}, which in term implies that
  $$
    \int^{i\infty}_{i\epsilon} z^m S_2(z)\,dz
  $$
  converges absolutely for any given $\epsilon>0$. It follows that
  \begin{equation*}
  \begin{split}
    \int_0^\infty z^m S_2(z)\,dz
   &=2\lim_{\epsilon\to 0}
    \int^{i\infty}_{i\epsilon}\sum_{c\in\ZZ,c\neq 0}\frac{dz}{(cNz+1)^{n+1}}\\
   &=2\lim_{\epsilon\to 0}
    \sum_{c\in\ZZ,c\neq 0} \int^{i\infty}_{i\epsilon}
    \frac{dz}{(cNz+1)^{n+1}}
   =\lim_{\epsilon\to 0}\frac 2{Nn}\sum_{c\in\ZZ,c\neq 0}
    \frac1{c(icN\epsilon+1)^n} \\
  &=\lim_{\epsilon\to 0}\frac{2\epsilon}{Nn}\sum_{c=1}^\infty
    \frac{(1-icN\epsilon)^n-(1+icN\epsilon)^n}
    {(c\epsilon)(1+c^2N^2\epsilon^2)^n}.
  \end{split}
  \end{equation*}
  Now the last expression is a Riemann sum of the integral
  \begin{equation} \label{eq lost track}
    2\int^\infty_0\frac{(1-iNx)^n-(1+iNx)^n}{Nnx(1+N^2x^2)^n}\,dx
 =\int^\infty_{-\infty}
  \frac{(1-iNx)^n-(1+iNx)^n}{Nnx(1+N^2x^2)^n}\,dx,
  \end{equation}
  which, by moving the line of integration to $\mathrm{Im}\,x=\infty$ and
  counting the residue at $i/N$, is shown to be equal to $-2\pi i/(Nn)$.
  Thus, we conclude that
  \begin{equation}\label{eqn2.4}
    \int_0^{i\infty}z^m S_2(z)\,dz
   =\begin{cases}
    0, &\text{if }m>\tn+1, \\
    -2\pi i/(Nn), &\text{if }m=\tn+1. \end{cases}
  \end{equation}

  We now consider the contribution from $S_3(z)$. For this it is
  necessary to distinguish the two cases $m<w$ and $m=w$. For the
  former case $m<w$, again, we first show that the integral
  $$
    \int_0^{i\infty}z^m S_3(z)\,dz=\int^{i\infty}_0 z^m\sum_{bc\neq 0}
    \frac{dz}{(az+b)^{\tn+1}(cz+d)^{n+1}}
  $$
  converges absolutely. For $|z|\ll 1$, write
  $$
    S_3(z)=\sum_{b,d}\frac1{b^{\tn+1}d^{n+1}}\sum_h
    \left(\frac{cz}d+\frac z{bd}+hNz+1\right)^{-(\tn+1)}
    \left(\frac cdz+hNz+1\right)^{-(n+1)},
  $$
  where in the inner sum, for a given pair of $b$ and $d$, the number
  $c$ is an integer satisfying $N|c$, $|c|\le dN$, and $bc\equiv
  -1\mathrm{\ mod\ }d$, and $h$ runs over all integers satisfying
  $c/d+hN\neq 0$. (In particular, when $d\neq\pm 1$, $h$ runs over all
  integers.) Now
  \begin{equation*}
  \begin{split}
   &\left|\frac{cz}d+\frac z{bd}+hNz+1\right|^{-(\tn+1)}
    \left|\frac cdz+hNz+1\right|^{-(n+1)} \\
   &\qquad\ll\begin{cases}
      1, &\text{if }|hN+c/d|\ll 1/|z|, \\
      |(hN+c/d)z|^{-w-2}, &\text{if }|hN+c/d|\gg 1/|z|.
    \end{cases}
  \end{split}
  \end{equation*}
  It follows that
  \begin{equation} \label{eq inserted 2.5}
    \sum_{\substack{\left[
      \begin{smallmatrix}a&b\\ c&d\end{smallmatrix}\right]
      \in\iGamma_0(N)},bc\neq 0}
    \frac1{|az+b|^{\tn+1}|cz+d|^{n+1}}\ll\frac 1{|z|}
  \end{equation}
  for $|z|\ll 1$. Furthermore, for $|z|\gg 1$, we write
  \begin{equation*}
  \begin{split}
    &\sum_{\substack{\left[
       \begin{smallmatrix}a&b\\ c&d\end{smallmatrix}\right]
       \in\iGamma_0(N)},bc\neq 0}
     \frac1{|az+b|^{\tn+1}|cz+d|^{n+1}} \\
    &\qquad=\frac1{N^{n+1}|z|^{w+2}}\sum_{\substack{\left[
     \begin{smallmatrix}d&-c/N\\ -bN&a\end{smallmatrix}\right]
     \in\iGamma_0(N)},bc\neq 0}
     \frac1{|a+bN/(Nz)|^{\tn+1}|c/N+d/(Nz)|^{n+1}}.
  \end{split}
  \end{equation*}
  Then by the same argument as before, we find that the sum is bounded
  by $|z|^{-w-1}$. It follows that
  $$
    \int^{i\infty}_0 z^mS_3(z)\,dz
  $$
  converges absolutely when $m<w$ and we may change the order of
  integration and summation freely. Now we have
  \begin{equation*}
  \begin{split}
    S_3(-z)&=\sum_{\substack{\left[
      \begin{smallmatrix}a&b\\ c&d\end{smallmatrix}\right]
      \in\iGamma_0(N)},bc\neq 0}\frac1{(-az+b)^{\tn+1}(-cz+d)^{n+1}} \\
    &=(-1)^{\tn+1}\sum_{\substack{\left[
      \begin{smallmatrix}a&b\\ c&d\end{smallmatrix}\right]
      \in\iGamma_0(N)},bc\neq 0}\frac1{(az-b)^{\tn+1}(-cz+d)^{n+1}}
     =(-1)^{\tn+1}S_3(z)
  \end{split}
  \end{equation*}
  and thus
  \begin{equation} \label{eq inserted 2.6}
    \int^{i\infty}_0 z^mS_3(z)\,dz=\frac12\int_{-i\infty}^{i\infty}
    z^m S_3(z)\,dz.
  \end{equation}
  Furthermore, for each $\left[\begin{smallmatrix}a&b\\
    c&d\end{smallmatrix}\right]\in\iGamma_0(N)$ with $bc\neq 0$ we have
  $abcd>0$ and consequently the two poles $-b/a$ and $-d/c$ of
  $z^m/(az+b)^{\tn+1}(cz+d)^{n+1}$ lie on the same side of the imaginary
  axis. This implies that
  \begin{equation} \label{eq inserted 2.7}
    \int^{i\infty}_{-i\infty}\frac{z^m}{(az+b)^{\tn+1}(cz+d)^{n+1}}\,dz
   =0
  \end{equation}
  for all $\left[\begin{smallmatrix}a&b\\
    c&d\end{smallmatrix}\right]\in\iGamma_0(N)$ with $bc\neq 0$ and we
  have
  \begin{equation}\label{eqn2.5}
    \int^{i\infty}_0 z^mS_3(z)\,dz=0
  \end{equation}
  when $m<w$.

  We now consider the contribution of $S_3$ for the case $m=w$. We
  write, by \eqref{eq inserted 2.6},
  \begin{equation*}
  \begin{split}
    \int^{i\infty}_0 z^wS_3(z)\,dz
  &=\frac12\lim_{\epsilon\to 0}\int^{i/\epsilon}_{-i/\epsilon}
     z^wS_3(z)\,dz \\
  &=\frac12\lim_{\epsilon\to 0}\sum_{\substack{\left[
      \begin{smallmatrix}a&b\\ c&d\end{smallmatrix}\right]
      \in\iGamma_0(N)},bc\neq 0}\int^{i/\epsilon}_{-i/\epsilon}
    \frac{z^w\,dz}{(az+b)^{\tn+1}(cz+d)^{n+1}}.
  \end{split}
  \end{equation*}
  Here the change of order of summation and integration is justified
  by the estimate \eqref{eq inserted 2.5}. Using
  \eqref{eq inserted 2.7} we reduce the last expression to
  $$
    -\frac12\lim_{\epsilon\to 0}\sum_{\substack{\left[
     \begin{smallmatrix}a&b\\ c&d\end{smallmatrix}\right]
     \in\iGamma_0(N)},bc\neq 0}\left(\int^{i\infty}_{i/\epsilon}
    +\int^{-i/\epsilon}_{-i\infty}\right)
     \frac{z^w\,dz}{(az+b)^{\tn+1}(cz+d)^{n+1}}.
  $$
  We then make a change of variable $z=i/(\epsilon t)$ and obtain
  $$
    \int^{i\infty}_0 z^wS_3(z)\,dz
  =-\lim_{\epsilon\to 0}\sum_{\substack{\left[
     \begin{smallmatrix}a&b\\ c&d\end{smallmatrix}\right]
     \in\iGamma_0(N)},bc\neq 0}
    \frac{\epsilon}{2i}\int^1_{-1}
    \frac{dt}{(a-i\epsilon tb)^{\tn+1}(c-i\epsilon td)^{n+1}}.
  $$
  We rewrite the sum as
  $$
    \sum_{\substack{\left[
      \begin{smallmatrix}a&b\\ c&d\end{smallmatrix}\right]
      \in\iGamma_0(N)},c\neq 0}
   -\sum_{\substack{\left[
      \begin{smallmatrix}1&0\\ c&1\end{smallmatrix}\right]
      \in\iGamma_0(N)},c\neq 0}
   -\sum_{\substack{\left[
      \begin{smallmatrix}-1&0\\ c&-1\end{smallmatrix}\right]
      \in\iGamma_0(N)},c\neq 0}.
  $$
  The latter two sums are bounded above by
  $$
    \ll\epsilon\sum_{c=1}^\infty\int^1_{-1}
    \frac{dt}{|c-i\epsilon t|^{n+1}}\ll\epsilon
    \sum_{c=1}^\infty\frac1{c^{n+1}}\ll\epsilon.
  $$
  Thus, as $\epsilon$ tends to $0$, the latter two sums vanish.
  Consequently we have
  \begin{equation*}
  \begin{split}
    \int_0^{i\infty} z^wS_3(z)\,dz
 &=-\lim_{\epsilon\to 0}\sum_{\substack{\left[
      \begin{smallmatrix}a&b\\ c&d\end{smallmatrix}\right]
      \in\iGamma_0(N)},\,c\neq 0}
    \frac{\epsilon}{2i}\int^1_{-1}
    \frac{dt}{(a-i\epsilon tb)^{\tn+1}(c-i\epsilon td)^{n+1}} \\
 &=-\lim_{\epsilon\to 0}\sum_{(a,c)=1,\,N|c,\,c\neq 0}
    \frac1{a^{\tn+1}c^{n+1}} \\
 &\qquad\times
    \sum_{h\in\mathbb Z+b/a}\frac{\epsilon}{2i}
    \int^1_{-1}\frac{dt}{(1-i\epsilon th)^{\tn+1}
    (1-i\epsilon t(h+1/ac))^{n+1}},
  \end{split}
  \end{equation*}
  where for a given pair $(a,c)$, $b$ denotes an integer
  satisfying $bc\equiv-1\mod a$. Now the inner sum is a Riemann sum of
  $$
    \frac1{2i}\int^\infty_{-\infty}\int^1_{-1}
    \frac{dt}{(1-ixt)^{w+2}}\,dx.
  $$
  The evaluation of this integral is similar to that of
  \eqref{eq lost track} and we find that it equals
  $$
   -\frac1{2(w+1)}\int^\infty_{-\infty}
    \frac1x\left(\frac1{(1-ix)^{w+1}}-\frac1{(1+ix)^{w+1}}\right)\,dx
  =-\frac{\pi i}{w+1}.
  $$
  Therefore, we have
  $$
    \int^{i\infty}_0 z^wS_3(z)\,dz=\frac{\pi i}{w+1}
    \sum_{(a,c)=1,\,N|c,\,c\neq 0}\frac1{a^{\tn+1}c^{n+1}}.
  $$
  Finally, we find
  \begin{equation*}
  \begin{split}
    \sum_{(a,c)=1,\,N|c,\,c\neq 0}\frac1{a^{\tn+1}c^{n+1}}
    &=4\sum_{a=1,\,(a,N)=1}^\infty\frac1{a^{\tn+1}}
      \sum_{c=1,\,(a,c)=1}^\infty\frac1{(cN)^{n+1}} \\
    &=\frac4{N^{n+1}}\sum_{a=1,\,(a,N)=1}^\infty\frac1{a^{\tn+1}}
      \sum_{d|a}\mu(d)\sum_{c=1}^\infty\frac1{(cd)^{n+1}} \\
    &=\frac{4\zeta(n+1)}{N^{n+1}}\sum_{d=1,\,(d,N)=1}^\infty
      \frac{\mu(d)}{d^{n+1}}\sum_{a=1,\,(a,N)=1}^\infty\frac1{(ad)^{\tn+1}}\\
    &=\frac{4\zeta(n+1)\zeta(\tn+1)}{\zeta(w+2)N^{n+1}}\prod_{p|N}
      \frac{1-p^{-(\tn+1)}}{1-p^{-(w+2)}},
  \end{split}
  \end{equation*}
  where the product is taking over all prime divisors $p$ of $N$. With
  \eqref{eqn2.1} the last expression becomes
  $$
    -\frac2{N^{n+1}}\binom{w+2}{n+1}\frac{B_{n+1}B_{\tn+1}}{B_{w+2}}
     \prod_{p|N}\frac{1-p^{-(\tn+1)}}{1-p^{-(w+2)}}.
  $$
  In summary, we see that, for integers $m$ and $n$, $0\le m\le w$,
  $0<n<w$, with opposite parity and $m+n>w$,
  \begin{equation} \label{eq inserted 2.10}
    \int^{i\infty}_0 z^mS_3(z)\,dz
   =\begin{cases} 0, &\text{if }m<w, \\
   -\displaystyle\binom{w+2}{n+1}
    \frac{2\pi iB_{n+1}B_{\tn+1}}{(w+1)N^{n+1}B_{w+2}}
    \prod_{p|N}\frac{1-p^{-(\tn+1)}}{1-p^{-(w+2)}}, &\text{if }m=w.
    \end{cases}
  \end{equation}
  Combining \eqref{eqn2.2}, \eqref{eqn2.4} and \eqref{eq inserted 2.10},
  we conclude that
  \begin{equation*}
  \begin{split}
    c_{w,n}r_m(R_{\iGamma_0(N),w,n})&=\frac{2\pi i}{m+1}\binom{m+1}\tn
    B_{m-\tn+1}-\frac{2\pi i}{Nn}\delta_{m,\tn+1} \\
    &\qquad-\delta_{m,w}\binom{w+2}{n+1}
     \frac{2\pi iB_{n+1}B_{\tn+1}}{(w+1)N^{n+1}B_{w+2}}
     \prod_{p|N}\frac{1-p^{-(\tn+1)}}{1-p^{-(w+2)}}
  \end{split}
  \end{equation*}
  for the cases $m+n>w$.

  For integers $m$ and $n$ satisfying $m\not\equiv n\mod 2$ and $m<\tn$,
  or equivalently, $m+n<w$, we make use of Lemma \ref{lem2.2}. We have
  $$
    r_m(R_n)=(-N)^{\tn-m}r_{\tm}(R_\tn).
  $$
  Now $\tm+\tn>w$ and we can apply the formula just obtained with $m$ and
  $n$ replaced by $\tm$ and $\tn$, respectively. We find that
  \begin{equation*}
  \begin{split}
    c_{w,n}r_m(R_n)&=(-N)^{\tn-m}
    \Bigg(\frac{2\pi i}{\tm+1}\binom{\tm+1}nB_{\tm-n+1}
   -\frac{2\pi i}{N\tn}\delta_{\tm,n+1} \\
   &\qquad\qquad-\delta_{m,0}\binom{w+2}{n+1}
    \frac{2\pi iB_{n+1}B_{\tn+1}}{(w+1)N^{\tn+1}B_{w+2}}
    \prod_{p|N}\frac{1-p^{-(n+1)}}{1-p^{-(w+2)}}\Bigg).
  \end{split}
  \end{equation*}
  This proves the the proposition for the cases $0<n<w$.
  Now for the cases $n=0$ and $n=w$, by Definition \ref{defn1.1}, we have
  $$
    r_m(R_0)=2^{-1}(2i)^{w+1}(R_0,R_m)=2^{-1}(2i)^{w+1}\overline{(R_m,R_0)}
   =-\overline{r_0(R_m)}
  $$
  and also
  $$
    r_m(R_w)=-\overline{r_w(R_m)}
  $$
  At this point we insert the formulas just obtained into the equations
  above. After simplification, we arrive at the claimed formulas. This
  completes the proof.
  \end{proof}

  We now give a proof of Theorem \ref{thm1.1}.

  \begin{proof}[Proof of Theorem \ref{thm1.1}]
  By Proposition \ref{prop2.3}, we have
  \begin{equation*}
  \begin{split}
    r^-(R_{\iGamma_0(N),w,n})(X)
    &=-\sum_{\substack{0<m<w \\ \text{$m$\ odd}}}\binom wm X^m r_{\tm}
     (R_{\iGamma_0(N),w,n}) \\
    &=-\binom wn^{-1}\sum_{\substack{0<m<n \\ \text{$m$\ odd}}}
      \binom wm\binom{\tm+1}\tn\frac{X^m}{\tm+1}
      B_{\tm-\tn+1} \\
    &\qquad+\binom wn^{-1}\sum_{\substack{n<m<w \\ \text{$m$\ odd}}}
     \binom wm\binom{m+1}n\frac{N^{\tn-\tm}X^m}{m+1}B_{m-n+1} \\
    &\qquad+\binom wn^{-1}\binom w{n-1}\frac{X^{n-1}}{Nn}
    -\binom wn^{-1}\binom w{n+1}\frac{X^{n+1}}\tn.
  \end{split}
  \end{equation*}
  After a straightforward simplification, the first sum in the last
  expression becomes
  $$
    \frac1{n+1}\sum_{\substack{0<m<n \\ \text{$m$\ odd}}}
    \binom{n+1}mX^mB_{n-m+1}
   =\frac1{n+1}\left(B_{n+1}^0(X)-X^{n+1}\right),
  $$
  while the second sum is reduced to
  \begin{equation*}
  \begin{split}
   &\frac1{\tn+1}\sum_{\substack{n<m<w \\ \text{$m$\ odd}}}\binom{\tn+1}{\tm}
    N^{\tn-\tm}X^{w-\tm}B_{\tn-\tm+1} \\
   &\qquad\qquad=\frac{N^{\tn}X^w}{\tn+1}\left\{B_{\tn+1}^0
    \left(\frac1{NX}\right)-\frac1{(NX)^{\tn+1}}\right\}.
  \end{split}
  \end{equation*}
  Altogether, we see that
  $$
    r^-(R_{\iGamma_0(N),w,n})(X)=
    \frac{N^\tn X^w}{\tn+1}B_{\tn+1}^0(\frac1{NX})-\frac1{n+1}B_{n+1}^0(X),
  $$
  as claimed in the statement.

  The proof of the statement about $r^+(R_{\iGamma_0(N),w,n})$ for odd
  $n$ is almost the same, except that there are two extra terms
  $$
    c_{w,n}^{-1}\binom{w+2}{n+1}
    \frac{2\pi iB_{n+1}B_{\tn+1}}{(w+1)B_{w+2}}
    \left(\frac1{N^{n+1}}\prod_{p|N}\frac{1-p^{-(\tn+1)}}{1-p^{-(w+2)}}
   -\frac{X^w}{N}
    \prod_{p|N}\frac{1-p^{-(n+1)}}{1-p^{-(w+2)}}\right).
  $$
  We shall omit the details here.
  \end{proof}
\end{section}

\begin{section}{Proof of Lemma \ref{thm1.5} and Theorem \ref{thm1.6}}
  \label{sect3}
  Let $m$ be a positive integer. Recall that the congruence subgroup
  $\iGamma_0(N)$ acts on the set
  \begin{equation*}
    H_{N,m}:=\left\{
      \begin{bmatrix}a&b\\c&d\end{bmatrix}
      \Big|\ a,b,c,d\in\ZZ;\ ad-bc=m;\
          c\equiv 0 \pmod {N};\ \gcd(a,N)=1
      \right\}
  \end{equation*}
  by matrix multiplication on the left. The standard set of coset
  representatives is given by
  \begin{equation*}
    M_{N,m}:=\left\{
       \begin{bmatrix}a&b\\0&d\end{bmatrix}
       \ \Big|\ \ a,d\in\ZZ^+;\ ad=m;\
              \gcd(a,N)=1;\
              b=0,\ldots,d-1;\
       \right\}.
  \end{equation*}
  In particular, we have the following lemma.
  \begin{lem}\label{lem4.1} One has
  \begin{equation*}
    \iGamma_0(N) M_{N,m}=H_{N,m}.
  \end{equation*}
  \end{lem}
  With this lemma, it is straightforward to prove Lemma \ref{thm1.5}.
  \begin{proof}[Proof of Lemma \ref{thm1.5}]
  Let $f\in S_{w+2}(\iGamma_0(N))$ be a cusp form. Recall that the
  Hecke operator $T_m$ is defined by
  $$
    T_mf(z)=m^{w/2}\sum_{\alpha\in M_{N,m}}(f|_\alpha)(z).
  $$
  Thus, from Lemma \ref{lem2.1}, we have
  \begin{align*}
    T_mR_{\iGamma_0(N),w,n}(z)
    &=c_{w,n}^{-1}
    T_m\sum_{\substack{\left[\begin{smallmatrix}a&b\\
      c&d\end{smallmatrix}\right]\in \iGamma_0(N)}}
           \frac{1}{(az+b)^{\tn+1}(cz+d)^{n+1}} \\
    &=c_{w,n}^{-1}T_m\sum_{\gamma \in \iGamma_0(N)}
           z^{-\tn-1}|_\gamma \\
    &=m^{w/2}c_{w,n}^{-1}\sum_{\alpha \in
      M_{N,m}}\sum_{\gamma \in \iGamma_0(N)}
           z^{-\tn-1}|_\gamma|_\alpha \\
    &=m^{w/2}c_{w,n}^{-1}\sum_{\gamma' \in H_{N,m}}
           z^{-\tn-1}|_{\gamma'} \\
    &=m^{w+1}c_{w,n}^{-1}\sum_{\substack
             {\left[\begin{smallmatrix}a'&b'\\c'&d'\end{smallmatrix}\right]
                      \in H_{N,m}}}
           \frac{1}{(a'z+b')^{\tn+1}(c'z+d')^{n+1}} \\
    &=R_{\iGamma_0(N),w,n}^{m}(z).
  \end{align*}
  Here we used the identity
  \begin{equation*}
    H_{N,m}=\iGamma_0(N)M_{N,m}
  \end{equation*}
  in Lemma \ref{lem4.1}. This completes the proof.
  \end{proof}

  \begin{proof}[Proof of Theorem \ref{thm1.6}]
  The proofs of Proposition \ref{prop2.3} and Theorem \ref{thm1.1}
  are valid for this theorem with a few modifications.
  Let $n$ be an even integer and write
  $$
    R_n^m(z)=R_{\iGamma_0(N),w,n}^m
   =\frac{m^{w+1}}{c_{w,n}}\sum_{\substack{
    \left[\begin{smallmatrix}a&b\\c&d\end{smallmatrix}\right]\in
    H_{N,m}}}\frac1{(az+b)^{\tn+1}(cz+d)^{n+1}}.
  $$
  We are required to evaluate the integrals
  $$
    \int_0^{i\infty}z^jR_n^m(z)\,dz
  $$
  for odd integers $j$.

  For $j$ with $j>\tn$, we partition the sum $R_n^m(z)$
  as
  $$
    m^{-w-1}c_{w,n}R_n^m(z)=\sum_{c=0}+\sum_{b=0,c\neq0}
   +\sum_{abcd>0}+\sum_{d=0}+\sum_{abcd<0}
   =S_1+S_2+S_3+S_4+S_5.
  $$
  The evaluation of $\int^{i\infty}_0S_i(z)z^j\,dz$, $i=1,2,3$,
  is almost the same as that of the corresponding integrals in
  Proposition \ref{prop2.3}. For instance, if we replace $z$ by $az$
  in \eqref{S1}, multiply by $d^{-(n+1)}$ and sum over all pairs of
  $(a,d)$ satisfying $ad=m$, then we get the $S_1(z)$ function here.
  Thus, it is easy to see that
  \begin{equation} \label{tt S1}
    \int^{i\infty}_0z^jS_1(z)\,dz
   =\frac{2\pi i}{j+1}\binom{j+1}{\tn}B_{j-\tn+1}
    \sum_{\substack{ad=m,\ a>0\\ \gcd(a,N)=1}}\frac1{a^{j+1}d^{n+1}}
  \end{equation}
  for $j>\tn$. Likewise, for $S_2(z)$ we have
  \begin{equation} \label{tt S2}
    \int^{i\infty}_0z^jS_2(z)\,dz
   =\begin{cases} 0, &\text{if }j>\tn+1, \\
    \displaystyle-\frac{2\pi i}{Nn}\sum_{
    \substack{ad=m,\ a>0\\ \gcd(a,N)=1}}\frac1{a^{\tn+1}d^n},
    &\text{if }j=\tn+1.
    \end{cases}
  \end{equation}
  Furthermore, the argument in Proposition \ref{prop2.3} shows that
  \begin{equation} \label{tt S3}
    \int^{i\infty}_0z^jS_3(z)\,dz=0
  \end{equation}
  since $j$ is assumed to be odd and can not be $0$ or $w$.

  The sum $S_4$ is actually empty unless $N|m$. When this occurs, we
  have
  \begin{equation*}
  \begin{split}
    S_4(z)&=\sum_{\substack{\left[\begin{smallmatrix}
    a&b\\ cN& 0\end{smallmatrix}\right]\in H_{N,m}}}
    \frac1{(az+b)^{\tn+1}(cNz)^{n+1}} \\
    &=2\sum_{\substack{c>0,\\-bcN=m}}\frac1{(cNz)^{n+1}}
    \sum_{\substack{h=1,\\ \gcd(h,N)=1}}^N
    \sum_{a\in\ZZ}\frac1{((aN+h)z+b)^{\tn+1}}.
  \end{split}
  \end{equation*}
  The innermost sum is equal to
  \begin{equation*}
    \frac1{(Nz)^{\tn+1}}\sum_{a\in\ZZ}\frac1{(a+h/N+b/Nz)^{\tn+1}}
   =\frac{(-2\pi i)^{\tn+1}}{(Nz)^{\tn+1}\tn!}
    \sum_{r=1}^\infty r^\tn e^{2\pi ir(h/N+b/Nz)}.
  \end{equation*}
  Set
  \begin{equation*}
  \begin{split}
    f(r)&=\sum_{\substack{h=1\\ \gcd(h,N)=1}}^N e^{2\pi irh/N}
    =\sum_{d|N}\mu(d)\sum_{h=1}^{N/d}e^{2\pi irhd/N} \\
   &=N\sum_{d|N,\ N/d|r}\frac{\mu(d)}d=\sum_{d|\gcd(N,r)}
     \mu(N/d)d.
  \end{split}
  \end{equation*}
  Then we have
  \begin{equation*}
  \begin{split}
    \int^{i\infty}_0z^jS_4(z)\,dz
  &=\frac{2(-2\pi i)^{\tn+1}}{N^{w+2}\tn!}\int^{i\infty}_0
    z^{-(\tj+2)}\sum_{\substack{c>0\\ -bcN=m}}\frac1{c^{n+1}}
    \sum_{r=1}^\infty f(r)r^\tn e^{2\pi irb/Nz}\,dz \\
  &=\frac{2(-2\pi i)^{\tn+1}}{N^{w+2}\tn!}\int^{i\infty}_0
    z^{\tj}\sum_{\substack{c>0\\ -bcN=m}}\frac1{c^{n+1}}
    \sum_{r=1}^\infty f(r)r^\tn e^{2\pi ir|b|z/N}\,dz,
  \end{split}
  \end{equation*}
  where $\tj=w-j$. When $n>j$, we may integrate term by term, and the
  result is
  \begin{equation*}
  \begin{split}
    \int^{i\infty}_0z^jS_4(z)\,dz
   =\frac{2(-2\pi i)^{j-n}}{N^{j+1}}\frac{\tj!}{\tn!}
    \sum_{\substack{c>0 \\ -bcN=m}}\frac1{c^{n+1}|b|^{\tj+1}}
    \sum_{r=1}^\infty\frac{f(r)}{r^{n-j+1}}.
  \end{split}
  \end{equation*}
  Now we have
  $$
    \sum_{r=1}^\infty\frac{f(r)}{r^{n-j+1}}
   =\sum_{d|N}\mu(N/d)d\sum_{r=1}^\infty\frac1{(dr)^{n-j+1}}
   =\zeta(n-j+1)\sum_{d|N}\mu(N/d)d^{j-n}.
  $$
  Simplifying the result, we get
  \begin{equation*}
  \begin{split}
    \int^{i\infty}_0z^jS_4(z)\,dz
  &=\frac{2\pi iB_{n-j+1}\tj!}{N^{j+1}(n-j+1)!\tn!}
    \sum_{\substack{c>0\\ -bcN=m}}\frac1{c^{n+1}|b|^{\tj+1}}
    \sum_{d|N}\mu(N/d)d^{j-n} \\
  &=\frac{2\pi iN^{w-2j}B_{n-j+1}}{(\tj+1)m^{\tj+1}}\binom{\tj+1}\tn
    \sum_{c|(m/N)}c^{\tj-n}\sum_{d|N}\mu(N/d)d^{j-n}
  \end{split}
  \end{equation*}
  for $n>j$.

  For $n<j$, we write
  $$
    \int^{i\infty}_0z^jS_4(z)\,dz
   =i^{j+1}\int^\infty_0t^jS_4(it)\,dt
  $$
  and observe that the integral
  $$
    \int^\infty_0 t^sS_4(it)\,dt
  $$
  defines an analytic function for $\mathrm{Re}\,s<n$. Thus, by the
  uniqueness of analytic continuation, the formula
  \begin{equation*}
  \begin{split}
    \int^{i\infty}_0z^jS_4(z)\,dz
   =\frac{2(-2\pi i)^{j-n}}{N^{j+1}}\frac{\tj!}{\tn!}\zeta(n-j+1)
    \sum_{\substack{c>0 \\ -bcN=m}}\frac1{c^{n+1}|b|^{\tj+1}}
    \sum_{d|N}\mu(N/d)d^{j-n}.
  \end{split}
  \end{equation*}
  remains valid for $n<j$. In summary, we have shown that
  \begin{equation} \label{tt S4}
  \begin{split}
    \int^{i\infty}_0z^jS_4(z)\,dz
   =\frac{2\pi iN^{w-2j}B_{n-j+1}}{(\tj+1)m^{\tj+1}}
    \binom{\tj+1}\tn\sum_{c|(m/N)}c^{\tj-n}
    \sum_{d|N}\mu(N/d)d^{j-n},
  \end{split}
  \end{equation}
  where $B_\ell$ is understood to be $0$ if $\ell<0$.

  For $S_5$, we have
  $$
   \sum_{abcd<0}\int^{i\infty}_0\frac{z^j}{(az+b)^{\tn+1}(cz+d)^{n+1}}
   \,dz=\frac12\sum_{abcd<0}\int^{i\infty}_{-i\infty}
   \frac{z^j}{(az+b)^{\tn+1}(cz+d)^{n+1}}\,dz.
  $$
  To evaluate the integrals, it turns out to be easier to work with
  the polynomial
  $$
    F(X)=\int^{i\infty}_{-i\infty}\frac{(X-z)^w}
    {(az+b)^{\tn+1}(cz+d)^{n+1}}\,dz
  $$
  instead. Observe that
  $$
    \frac{d^\ell}{dX^\ell}F(X)\Big|_{X=-b/a}=0
  $$
  for $\ell=0,1,\ldots,n-1$, and
  $$
    \frac{d^\ell}{dX^\ell}F(X)\Big|_{X=-d/c}=0
  $$
  for $\ell=0,1\ldots,\tn-1$. This is because in these cases the
  integrands have only one pole at $-d/c$ for the former cases and at
  $-b/a$ for the latter cases. This implies that
  \begin{equation} \label{thm1.6 1}
    F(X)=C(aX+b)^n(cX+d)^\tn
  \end{equation}
  for some constant $C$. To determine the constant $C$, we consider
  $$
    \frac{d^n}{dX^n}F(X)\Big|_{X=-b/a}.
  $$
  Using the integral definition of $F(X)$, we have
  \begin{equation*}
  \begin{split}
    \frac{d^n}{dX^n}F(X)\Big|_{X=-b/a}
   &=a^{-\tn}\frac{w!}{\tn!}\int^{i\infty}_{-i\infty}
    \frac{dz}{(az+b)(cz+d)^{n+1}} \\
   &=2\pi i\sgn(ab)a^{-\tn-1}\frac{w!}{\tn!}(-bc/a+d)^{-n-1} \\
   &=2\pi i\sgn(ab)m^{-n-1}a^{n-\tn}\frac{w!}{\tn!}.
  \end{split}
  \end{equation*}
  On the other hand, using \eqref{thm1.6 1} we have
  $$
    \frac{d^n}{dX^n}F(X)\Big|_{X=-b/a}=Cn!a^n(-bc/a+d)^{\tn}
   =Cm^\tn a^{n-\tn}n!.
  $$
  Comparing the two expressions, we see that
  $C=m^{-w-1}\sgn(ab)c_{w,n}$ and
  \begin{equation} \label{tt S5}
    \int^{i\infty}_0(X-z)^w S_5(z)\,dz
   =\frac12m^{-w-1}c_{w,n}\sum_{\substack{\left[
    \begin{smallmatrix}a&b\\ c&d\end{smallmatrix}\right]\in H_{N,m}\\
    abcd<0}}\sgn(ab)(aX+b)^n(cX+d)^\tn.
  \end{equation}
  Altogether, \eqref{tt S1}, \eqref{tt S2}, \eqref{tt S3},
  \eqref{tt S4}, and \eqref{tt S5} give the evaluation of
  $r_j(R_{\iGamma_0(N),w,n}^m)$ for $j>\tn$.

  For the cases $j<\tn$, we partition the sum
  $$
    m^{-w-1}c_{w,n}R_{\iGamma_0(N),w,n}^m(z)
  $$
  as
  $$
    \sum_{b=0}+\sum_{c=0,\ b\neq 0}+\sum_{abcd>0}
   +\sum_{d=0}+\sum_{abcd<0}=S_1+S_2+S_3+S_4+S_5
  $$
  instead. The terms $S_3$, $S_4$, and $S_5$ are the same as before.
  For $S_1$ and $S_2$ we first make a change of variable
  $z\mapsto-1/z$ in the integrals $\int^{\infty}_0z^jS_i(z)\,dz$
  first. Then the evaluation is done in the same way as above. We find
  $$
    \int^{i\infty}_0z^jS_1(z)\,dz=-\frac{2\pi iN^{\tn-j}}{\tj+1}
    \binom{\tj+1}n B_{\tj-n+1}\sum_{\substack{ad=m,\,a>0\\
    \gcd(a,N)=1}} \frac1{a^{\tn+1}d^{\tj+1}}
  $$
  and
  $$
    \int^{i\infty}_0z^jS_2(z)\,dz=
    \begin{cases}
    0, &\text{if }j<\tn-1, \\
    \displaystyle\frac{2\pi i}{\tn}
    \sum_{\substack{ad=m,\,a>0\\ \gcd(a,N)=1}}\frac1{a^\tn d^{n+1}},
    &\text{if }j=\tn+1.
    \end{cases}
  $$
  We then insert the all the estimates above into the definition of
  $r^-(R^m_{\iGamma_0(N),w,n})(X)$. After simplification we arrive at
  the claimed formula.
  \end{proof}
\end{section}

\begin{section}{Proof of Theorems \ref{thm1.2}, \ref{thm1.4}, and
  \ref{thm1.7}}
  \label{sect4}
  We first give a proof of Theorems \ref{thm1.2}.

  \begin{proof}[Proof of Theorems \ref{thm1.2}]
  Write $d=d_w$. To prove that
  $$
    \{R_{\iGamma_0(2),w,2i}\ |\ i=1,\ldots,d_w\}
  $$
  is a basis, it suffices to prove that the polynomials
  $$
    r^-(R_{\iGamma_0(2),w,2i})(X), \qquad
    i=1,2,\ldots, d
  $$
  are linearly independent over $\mathbb C$.

  By Theorem \ref{thm1.1}, for $i=1,2,\ldots,d$, the
  coefficients of $X^{w-2j+1}$, $j=1,2\ldots,d$, in
  $r^-(R_{\iGamma_0(2),w,2i})(X)$ is
  $$
    \frac{2^{w-2i-2j+1}}{w-2i+1}\binom{w-2i+1}{2j-1}B_{w-2i-2j+2}.
  $$
  Therefore, the proof of the theorem reduces to that of
  \begin{equation}\label{eqn3.1}
    \det_{1\le i,j\le d}\left[\frac{2^{w-2i-2j+1}}{w-2i+1}
    \binom{w-2i+1}{2j-1}B_{w-2i-2j+2}\right]\neq 0.
  \end{equation}
  This would follow immediately from formulas in
  Theorem A of Appendix A. Here we provide an alternative proof that
  does not require more advanced tools from the theory of continued
  fractions.

  Using \eqref{eqn2.1} we find
  $$
    \frac{2^{w-2i-2j+1}}{w-2i+1}\binom{w-2i+1}{2j-1}B_{w-2i-2j+2}
   =\frac{(-1)^{w/2-i-j}}{\pi^{w-2i-2j+2}}\frac{(w-2i)!}{(2j-1)!}
    \zeta(w-2i-2j+2)
  $$
  and consequently the determinant in \eqref{eqn3.1} is equal to
  $$
    C\det_{1\le i,j\le d}\left[\zeta(w-2i-2j+2)\right]
  $$
  for some non-zero number $C$. Now we have
  \begin{equation*}
  \begin{split}
    \det_{1\le i,j\le d}\left[\zeta(w-2i-2j+2)\right]
  &=\sum_{m_1,\ldots,m_d=1}^\infty
    \det_{1\le i,j\le d}\left[m_i^{-w+2i+2j-2}\right] \\
  &=\sum_{m_1,\ldots,m_d=1}^\infty
    \det_{1\le i,j\le d}\left[m_i^{2j-2}\right]\prod_{i=1}^d m_i^{-w+2i}.
  \end{split}
  \end{equation*}
  By the Vandermonde identity, the last expression is equal to
  $$
    \sum_{m_1,\ldots,m_d=1}^\infty\prod_{1\le i<j\le d}(m_i^2-m_j^2)
    \prod_{i=1}^d m_i^{-w+2i}.
  $$
  This sum can be alternatively written as
  $$
  \sum_{1\le m_1<\ldots<m_d}\prod_{1\le i<j\le d}(m_i^2-m_j^2)
  \sum_{\sigma\in S_d}\varepsilon(\sigma)\prod_{i=1}^d m_{\sigma(i)}^{-w+2i},
  $$
  where the inner sum runs over all elements in the permutation group
  $S_d$, and $\varepsilon(\sigma)$ denotes $+1$ or $-1$ according as
  the permutation $\sigma$ is even or odd. Observe that this inner sum
  is in fact the expansion of the determinant
  $$
    \det_{1\le i,j\le d}\left[m_j^{-w+2i}\right].
  $$
  Using the Vandermonde identity again, we evaluate the above
  determinant as
  $$
    \prod_{1\le i<j\le d}(m_i^2-m_j^2)\prod_{i=1}^d m_i^{2-w}.
  $$
  Therefore,
  $$
    \det_{1\le i,j\le d}\left[\zeta(w-2i-2j+2)\right]
   =\sum_{1\le m_1<\ldots<m_d}\prod_{1\le i<j\le d}(m_i^2-m_j^2)^2
    \prod_{i=1}^d m_i^{2-w}.
  $$
  Since the summands above are all strictly positive, this implies that
  \eqref{eqn3.1} holds. This proves that
  $\{R_{\iGamma_0(2),w,2i}\ |\ i=1,\ldots,d_w\}$ is a basis.

  The proofs of the other three cases are similar. For the set
  $\{R_{\iGamma_0(2),w,w-2i}\ |\ i=1,\ldots,d_w\}$, we find that the
  coefficients of $X^{2j-1}$ in $r^-(R_{\iGamma_0(2),w,w-2i})(X)$ is
  $$
   -\frac1{w-2i+1}\binom{w-2i+1}{2j-1}B_{w-2i-2j+2}.
  $$
  By the same argument as before, we see that $R_{\iGamma_0(2),w,w-2i}$,
  $i=1,\ldots,d_w$ are linearly independent. For the cases
  $\{R_{\iGamma_0(2),w,2i-1}\}$ and $\{R_{\iGamma_0(2),w,w-2i+1}\}$ we
  consider the coefficients of $X^{w-2j}$ and $X^{2j}$,
  $j=1,\ldots,d_w$, respectively, in $r^+(R_{\iGamma_0(2),w,2i-1})(X)$
  and $r^+(R_{\iGamma_0(2),w,w-2i+1})(X)$. The details are omitted.
\end{proof}

  \begin{proof}[Proof of Theorem \ref{thm1.4}] Write
  $k=w+2$. Let $S^\old_k(\iGamma_0(2))$ and $S^\new_k(\iGamma_0(2))$ be
  the spaces of oldforms and newforms of weight $k$ on $\iGamma_0(2)$,
  respectively. Let $d^\prime_k=\dim S_k(SL_2(\ZZ))$, and
  $f_1,\ldots,f_{d^\prime_k}$ be the normalized Hecke eigenforms
  spanning $S_k(SL_2(\ZZ))$. Then the functions
  $$
    g_i=f_i+f_i\big|_{\omega_2}, \quad
    g_{d^\prime_k+i}=f_i-f_i\big|_{\omega_2}, \quad
    i=1,\ldots,d^\prime_k,
  $$
  form an orthogonal basis for $S^\old_k(\iGamma_0(2))$. Let also
  $g_{2d^\prime_k+1},\ldots,g_{d_w}$ be the normalized Hecke
  eigenforms in $S^\new_k(\iGamma_0(2))$. Then the set
  $\{g_i\ |\ i=1,\ldots,d_w\}$ is an orthogonal basis for
  $S_k(\iGamma_0(2))$. It follows that
  $$
    E^0_{2j}E_{k-2j}^{i\infty}=\sum_{i=1}^{d_w}
    \frac{(E^0_{2j}E^{i\infty}_{k-2j},g_i)}{(g_i,g_i)}g_i,
  $$
  and to show that $\{E^0_{2j}E^{i\infty}_{k-2j}\}$ and
  $\{E^0_{k-2j}E^{i\infty}_{2j}\}$ are bases, it suffices to prove
  that
  $$
    \det_{\substack{1\le i\le d_w\\ 2\le j\le d_w+1}}\left[
    \frac{(E^0_{2j}E^{i\infty}_{k-2j},g_i)}{(g_i,g_i)}\right]\neq 0,
    \quad \det_{\substack{1\le i\le d_w\\ 2\le j\le d_w+1}}\left[
    \frac{(E^0_{k-2j}E^{i\infty}_{2j},g_i)}{(g_i,g_i)}\right]\neq 0.
  $$

  Now according to Proposition 2 of \cite{IK}, we have, for
  $i=1,\ldots,2d^\prime_k$,
  \begin{equation*}
  \begin{split}
    (E^0_{2j}E^{i\infty}_{k-2j},g_i)
  &=(E^0_{2j}E^{i\infty}_{k-2j},f_i\pm f_i\big|_{\omega_2}) \\
  &=c_j(1+2^{1-k}(1-\lambda_i))L(f_i,k-1)
    L(f_i\pm f_i\big|_{\omega_2},k-2j) \\
  &=c_j(1+2^{1-k}(1-\lambda_i))L(f_i,k-1)L(g_i,k-2j),
  \end{split}
  \end{equation*}
  and, for $i=2d^\prime_k+1,\ldots,d_w$,
  $$
    (E^0_{2j}E^{i\infty}_{k-2j},g_i)
   =c_j(\epsilon_i+2^{-k/2})L(g_i,k-1)L(g_i,k-2j)
  $$
  where $L(f,s)$ is the Hecke $L$-function associated with an
  eigenform $f$,
  $$
    c_j=\frac{(k-2)!}{(4\pi)^{k-1}}\cdot\frac{4j}{B_{2j}}\cdot
    \frac{2^j}{1-2^{2j}}\cdot\frac1{(1-2^{2j-k})\zeta(k-2j)},
  $$
  $\lambda_i=\lambda_{d^\prime_k+i}$ is the eigenvalue of $f_i$
  under the Hecke operator $T_2$, and $\epsilon_i$ is the eigenvalues
  of $g_i$ under the Atkin-Lehner involution $\omega_2$. For
  convenience, we have also set $f_i=f_{i-d^\prime_k}$ for
  $i=d^\prime_k+1,\ldots,2d^\prime_k$. It follows that
  $$
    \det_{\substack{1\le i\le d_w\\ 2\le j\le d_w+1}}\left[
    \frac{(E^0_{2j}E^{i\infty}_{k-2j},g_i)}{(g_i,g_i)}\right]
   =\det_{\substack{1\le i\le d_w\\ 2\le j\le d_w+1}}\left[
    L(g_i,k-2j)\right] \cdot\prod_{i=1}^{d_w}\frac{b_i}{(g_i,g_i)}\cdot
    \prod_{j=2}^{d_w+1}c_j
  $$
  and
  $$
    \det_{\substack{1\le i\le d_w\\ 2\le j\le d_w+1}}\left[
    \frac{(E^0_{k-2j}E^{i\infty}_{2j},g_i)}{(g_i,g_i)}\right]
   =\det_{\substack{1\le i\le d_w\\ 2\le j\le d_w+1}}\left[
    L(g_i,2j)\right] \cdot\prod_{i=1}^{d_w}\frac{b_i}{(g_i,g_i)}\cdot
    \prod_{j=2}^{d_w+1}c_{k/2-j},
  $$
  where
  $$
    b_i=\begin{cases}(1+2^{1-k}(1-\lambda_i))L(f_i,k-1),
    &\text{for }i=1,\ldots,2d^\prime_k, \\
    (\epsilon_i+2^{-k/2})L(g_i,k-1),
    &\text{for }i=2d^\prime_k+1,\ldots,d_w.
    \end{cases}
  $$
  Clearly, the numbers $c_j$ are never zero. Also, using the
  upper bound $|a_p|\le 2p^{(k-1)/2}$ for the $p$-th Fourier coefficients
  of a normalized Hecke eigenform $\sum a_nq^n$ on $SL_2(\ZZ)$, we
  know that $L(f_i,k-1)$, $L(g_i,k-1)$ and $(1+2^{1-k}(1-\lambda_i))$
  do not vanish.
  Then the proof of the theorem reduces to that of
  \begin{equation} \label{eq thm1.4 1}
    \det_{\substack{1\le i\le d_w\\ 2\le j\le d_w+1}}
    \left[L(g_i,k-2j)\right]\neq 0, \quad
    \det_{\substack{1\le i\le d_w\\ 2\le j\le d_w+1}}
    \left[L(g_i,2j)\right]\neq 0.
  \end{equation}

  On the other hand, consider the sets
  \begin{equation} \label{eq thm1.4 2}
    \{R_{\iGamma_0(2),w,2j+1}\,|\,j=1,\ldots,d_w\}, \quad
    \{R_{\iGamma_0(2),w,w-2j-1}\,|\,j=1,\ldots,d_w\}.
  \end{equation}
  Again, we have
  $$
    R_{\iGamma_0(2),w,n}=\sum_{i=1}^{d_w}
    \frac{(R_{\iGamma_0(2),w,n},g_i)}{(g_i,g_i)}g_i.
  $$
  \begin{equation*}
  \end{equation*}
  By Definition \ref{defn1.1},
  $$
    2^{-1}(2i)^{w+1}(g_i,R_{\iGamma_0(2),w,n})
    =\int^{i\infty}_0 g_i(z)z^n\,dz=\frac{n!}{(-2\pi i)^{n+1}}L(g_i,n+1).
  $$
  It follows that \eqref{eq thm1.4 1} holds if and only if the two sets
  in \eqref{eq thm1.4 2} are both bases for $S_{w+2}(\iGamma_0(2))$.

  To show that
  $$
    \{R_{\iGamma_0(2),w,2j+1}\ |\  j=1,\ldots,d_w\}
  $$
  is a basis, we consider the coefficients of $X^{w-2i}$ in
  $r^+(R_{\iGamma_0(2),w,2j+1})$. When $w$ is a multiple of $4$, the
  coefficients are
  $$
    \frac{2^{w-2i-2j-1}}{w-2j}\binom{w-2j}{2i}B_{w-2i-2j},
  $$
  and the proof proceeds as in Theorem \ref{thm1.2}, and the details
  are skipped. When $w=4d+2$ with $d_w=d$ for some positive integer
  $d$, the coefficients are
  $$
    a_{ij}=\frac{2^{w-2i-2j-1}}{w-2j}\binom{w-2j}{2i}B_{w-2i-2j}
   -\delta_{d,i}\delta_{d,j}\frac{B_0}{2d+2}.
  $$
  Then we have
  \begin{equation*}
  \begin{split}
    \det_{1\le i,j\le d}[a_{ij}]
  &=\det_{1\le i,j\le d}\left[
    \frac{2^{w-2i-2j-1}}{w-2j}\binom{w-2j}{2i}B_{w-2i-2j}\right] \\
  &\qquad\quad-\frac1{2d+2}\det_{1\le i,j\le d-1}\left[
    \frac{2^{w-2i-2j-1}}{w-2j}\binom{w-2j}{2i}B_{w-2i-2j}\right].
  \end{split}
  \end{equation*}
  By the first and the third identities of Theorem A in Appendix A,
  the first term is equal to
  \begin{equation*}
  \begin{split}
   &2^{d(w-1)}\left(\prod_{i=1}^d2^{-2i}\right)^2
    \prod_{j=1}^d(w-2j-1)!\prod_{i=1}^d\frac1{(2i)!}
    \det_{1\le i,j\le d}\left[\frac{B_{w-2i-2j}}{(w-2i-2j)!}\right] \\
   &\qquad=2^{-d}\prod_{i=1}^d\frac{(w-2i-1)!}{(2i)!}
    \prod_{i=1}^{2d-1}(2i+1)^{-(2d-i)},
  \end{split}
  \end{equation*}
  while the second term is
  \begin{equation*}
  \begin{split}
    &\frac{2^{(d-1)(w-1)}}{d+2}\left(\prod_{i=1}^{d-1}2^{-2i}\right)^2
     \prod_{j=1}^{d-1}(w-2j-1)!\prod_{i=1}^{d-1}\frac1{(2i)!}
     \det_{1\le i,j\le d-1}\left[\frac{B_{w-2i-2j}}{(w-2i-2j)!}\right] \\
    &\qquad=\frac{2^{1-d}}{2d+2}\prod_{i=1}^{d-1}\frac{(w-2i-1)!}{(2i)!}
     d(2d+1)\prod_{i=1}^{2d-1}(2i+1)^{-(2d-i)}.
  \end{split}
  \end{equation*}
  The ratio of these two terms is
  $$
    \frac{d+1}{d(2d+1)}\frac{(w-2d-1)!}{(2d)!}=\frac{d+1}d,
  $$
  which is never equal to $1$. From this we conclude that $\det_{1\le
  i,j\le d}[a_{ij}]\neq 0$ and thus that
  $\{E^0_{2j+2}E^{i\infty}_{w-2j}\ |\ j=1,\ldots,d_w\}$ is a basis.
  The proof of the assertion about
  $\{E^{i\infty}_{2j+2}E^0_{w-2j}\ |\ j=1,\ldots,d_w\}$ is similar, and
  is omitted here.
  \end{proof}

  \begin{proof}[Proof of Theorem \ref{thm1.7}] The proof uses a
  straightforward argument in elementary linear algebra, and is skipped
  here. For more details, see Theorem 2.9 of \cite{F5}.
  \end{proof}
\end{section}

\begin{section}{Acknowledgment}
  The authors would like to thank Professor C. Krattenthaler for
  his permission to include his proof of formulas for Hankel
  determinants of Bernoulli numbers in the paper. The authors would
  also like to thank Professors H. H. Chan and N. Yui for reading
  earlier versions of the paper and giving valuable comments.
  Finally, the authors would like to express their gratitude to the
  anonymous referee for thorough reading of the manuscript.

  Yifan Yang was supported by Grant 95-2115-M-009-005 of the National
  Science Council (NSC) of the Republic of China (Taiwan). Part of
  this work was done while he was visiting Tsuda College, Japan.
  He would like to thank Tsuda College for the enormous
  hospitality.
\end{section}
\medskip

\begin{center}
\bf{Appendix A: Evaluation of Hankel determinants of Bernoulli
  numbers}
\end{center}

\centerline{by {\sc Christian Krattenthaler}}
\medskip

In this appendix we will derive exact formulas for Hankel determinants
formed by Bernoulli numbers. We utilize the following theorem.
\medskip

\noindent{\bf Theorem CF.} {\it Let $(\mu_k)_{k\ge 0}$ be a sequence with
  generating function $\sum_{k=0}^\infty\mu_kx^k$ written in the form
  $$
    \sum_{k=0}^\infty\mu_kx^k
   =\cfrac{\mu_0}{1+\cfrac{a_1x}{1+\cfrac{a_2x}{1+\ddots}}}.
  $$
  Then we have
  \begin{align*}
  \det_{0\le i,j\le n-1}[\mu_{i+j}]&=\mu_0^n
  (a_1a_2)^{n-1}(a_3a_4)^{n-2}\ldots(a_{2n-5}a_{2n-4})^2
  (a_{2n-3}a_{2n-2}), \\
  \det_{0\le i,j\le n-1}[\mu_{i+j+1}]&=(-1)^n\mu_0^na_1^n
  (a_2a_3)^{n-1}(a_4a_5)^{n-2}\ldots(a_{2n-4}a_{2n-3})^2
  (a_{2n-2}a_{2n-1}),
  \end{align*}
  and
  \begin{equation*}
  \begin{split}
    \det_{0\le i,j\le n-1}[\mu_{i+j+2}]&=\mu_0^na_1^n
  (a_2a_3)^{n-1}(a_4a_5)^{n-2}\ldots(a_{2n-4}a_{2n-3})^2
  (a_{2n-2}a_{2n-1}) \\
  &\qquad\times\sum_{0\le i_1-1<i_2-2<\ldots<i_n-n\le n}
   a_{i_1}a_{i_2}\cdots a_{i_n}.
  \end{split}
  \end{equation*}
}
\medskip

Here the first two formulas are from Theorem 7.2 of \cite{JT}, and the
third is from Theorem 30 of \cite{Kr}.
\medskip

\noindent{\bf Theorem A.} {\it Let $B_n$ denote the $n$-th Bernoulli
  number. Then we have
  \begin{align*}
    \det_{0\le i,j\le n-1}\left[\frac{B_{2i+2j+2}}{(2i+2j+2)!}\right]
   &=4^{-n^2}\prod_{i=1}^{2n-1}(2i+1)^{-(2n-i)}, \\
    \det_{0\le i,j\le n-1}\left[\frac{B_{2i+2j+4}}{(2i+2j+4)!}\right]
   &=(-1)^n4^{-n^2-n}9^{-n}\prod_{i=1}^{2n-1}(2i+3)^{-(2n-i)}, \\
    \det_{0\le i,j\le n-1}\left[\frac{B_{2i+2j+6}}{(2i+2j+6)!}\right]
   &=4^{-n^2-2n}(n+1)(2n+3)\prod_{i=1}^{2n+1}(2i+1)^{-(2n+2-i)}.
  \end{align*}
}
\medskip
\begin{proof} In view of Theorem CF, we should determine the continued
  fraction expansion for the generating function for the numbers
  $B_{2n}/(2n)!$, $n\ge 1$. Writing $t=x/2$, we have
  \begin{equation*}
  \begin{split}
  \sum_{k\ge 0}\frac{B_{k+2}}{(k+2)!}x^k
  &=\frac1{x^2}\left(\frac x{e^x-1}-1+\frac x2\right)
   =\frac1{4t^2}\frac{te^t+te^{-t}-e^t+e^{-t}}{e^t-e^{-t}} \\
  &=\frac1{4t^2}\frac{\displaystyle\sum_{n\ge 1}(2n)
    \frac{t^{2n+1}}{(2n+1)!}}
    {\displaystyle\sum_{n\ge 0}\frac{t^{2n+1}}{(2n+1)!}}
   =\frac1{12}\frac{\displaystyle
   _0F_1\left(-;\frac52;\frac{t^2}4\right)}
    {\displaystyle _0F_1\left(-;\frac32;\frac{t^2}4\right)}.
  \end{split}
  \end{equation*}
  Equivalently, by using the fact that $B_n=0$ for odd $n\ge 3$,
  $$
    \sum_{k\ge 0}\frac{B_{2k+2}}{(2k+2)!}x^k
   =\frac1{12}\frac{\displaystyle _0F_1\left(-;\frac52;\frac x{16}\right)}
    {\displaystyle _0F_1\left(-;\frac32;\frac x{16}\right)}.
  $$

  Now, Gauss proved the following continued fraction expansion
  $$
    \frac{_2F_1(a,b+1;c+1;z)}{_2F_1(a,b;c;z)}=
    \cfrac1{1+\cfrac{a_1z}{1+\cfrac{a_2z}{1+\ddots}}},
  $$
  where
  $$
    a_{2n-1}=-\frac{(a+n-1)(c-b+n-1)}{(c+2n-2)(c+2n-1)}, \quad
    a_{2n}=-\frac{(b+n)(c-a+n)}{(c+2n-1)(c+2n)}.
  $$
  (See \cite[Theorem 6.1]{JT}.) Replacing $a$ by $b$, $c$ by $3/2$, and
  $z$ by $x^2/(16b^2)$ and subsequently letting $b$ tend to infinity,
  we obtain
  $$
    \frac{\displaystyle_0F_1\left(-;\frac52;\frac x{16}\right)}
    {\displaystyle_0F_1\left(-;\frac32;\frac x{16}\right)}
   =\cfrac1{1+\cfrac{a_1x}{1+\cfrac{a_2x}{1+\ddots}}},
  $$
  where
  $$
    a_{2n-1}=\frac1{4(4n-1)(4n+1)}, \quad
    a_{2n}=\frac1{4(4n+1)(4n+3)},
  $$
  that is,
  $$
    a_n=\frac1{4(2n+1)(2n+3)}.
  $$
  We now apply Theorem CF with $\mu_k=B_{2k+2}/(2k+2)!$. The first
  identity of the theorem yields
  \begin{equation*}
  \begin{split}
    \det_{0\le i,j\le n-1}\left[\frac{B_{2i+2j+2}}{(2i+2j+2)!}\right]
   &=12^{-n}\prod_{k=1}^{n-1}
     \left(16(4k-1)(4k+1)^2(4k+3)\right)^{-(n-k)} \\
   &=4^{-n^2}3^{-n}\left(3^{-(n-1)}5^{-(2n-2)}\ldots
     (4n-3)^2(4n-1)\right) \\
   &=4^{-n^2}\prod_{i=1}^{2n-1}(2i+1)^{-(2n-i)},
  \end{split}
  \end{equation*}
  and the second identity gives
  \begin{equation*}
  \begin{split}
    \det_{0\le i,j\le n-1}\left[\frac{B_{2i+2j+4}}{(2i+2j+4)!}\right]
   &=(-720)^{-n}\prod_{k=1}^{n-1}
     \left(16(4k+1)(4k+3)^2(4k+5)\right)^{-(n-k)} \\
   &=(-1)^n4^{-n^2-n}9^{-n}\prod_{i=1}^{2n-1}(2i+3)^{-(2n-i)}.
  \end{split}
  \end{equation*}
  For the third case, we have
  \begin{equation*}
  \begin{split}
    \det_{0\le i,j\le n-1}\left[\frac{B_{2i+2j+6}}{(2i+2j+6)!}\right]
  &=(-1)^n\det_{0\le i,j\le n-1}
    \left[\frac{B_{2i+2j+4}}{(2i+2j+4)!}\right]\\
  & \qquad\qquad
    \times\sum_{0\le i_1-1<i_2-2<\ldots<i_n-n\le n}a_{i_1}\ldots a_{i_n}.
  \end{split}
  \end{equation*}
  Now
  \begin{equation*}
  \begin{split}
    &\sum_{0\le i_1-1<i_2-2<\ldots<i_n-n\le n}a_{i_1}\ldots a_{i_n} \\
    &\quad=\sum_{0\le i_1-1<\ldots<i_n-n\le n}
     \frac1{(2i_1+1)(2i_1+3)(2i_2+1)(2i_2+3)\ldots(2i_n+1)(2i_n+3)}.
  \end{split}
  \end{equation*}
  The denominator of each summand is a product of $2n$ distinct odd
  integers among the $2n+1$ odd integers ranging from $3$ to $4n+3$.
  The only missing odd integer must be of the form $4k+3$,
  corresponding to the choice
  $(i_1,\ldots,i_n)=(1,3,\ldots,2k-1,2k+2,\ldots,2n-2,2n)$.
  Thus, we see that
  \begin{equation*}
  \begin{split}
    &\sum_{0\le i_1-1<\ldots<i_n-n\le n}a_{i_1}\ldots a_{i_n}
    =\frac1{4^n(4n+3)!!}\sum_{k=0}^n(4k+3)
    =\frac{(2n+3)(n+1)}{4^n(4n+3)!!},
  \end{split}
  \end{equation*}
  where $(4n+3)!!$ denotes the product of all odd integers less than
  or equal to $4n+3$. After simplification, we arrive at the third
  identity. This completes the proof.
\end{proof}

\begin{center}
\bf{Appendix B: Computing matrices representing the Hecke operators
on $S_{w+2}(\iGamma_0(2))$ and their characteristic polynomials}
\end{center}
\label{appendix}

We will demonstrate a
Mathematica\footnote{Mathematica is a trademark of Wolfram Research,~Inc.}
program which, for given $w$
and $m$, yields a matrix representing
the Hecke operator $T_m$ on $S_{w+2}(\iGamma_0(2))$,
and its characteristic polynomial.
The program is a straightforward application of Theorem \ref{thm1.7}.

\vspace{3mm}
\begin{verbatim}
w=10;(* w:even positive *)
m=6;(* m:positive *)
nn=2;(* nn:modulus, fixed *)
dw=Quotient[w+2,4]-1;
If[OddQ[w]||(dw<=0)||(m<=0)||(nn!=2),
    Print["w odd or dw<=0 or m<=0 or nn!=2"];Exit[]];
ber0[n_Integer,x_]:=
    Module[{b},
      b=BernoulliB[n,x]-n*BernoulliB[1]*x^(n-1);
      Return[b];
      ];
smwn=Array[a,{2,dw}];(* period polynomials *)
cosmwn=Array[d,{2,dw,w+1}];(* coefficients *)
For[j=1,j<=2,j++,
  If[j==1,em=1,em=m];
  For[i=1,i<=dw,i++,
    n=2*i; tn=w-n;
    fsm=0;
    For[a=-em+1,a<em,a++,
      For[b=-em+1,b<em,b++,
          For[c=-em+1,c<em,c++,
              For[d=-em+1,d<em,d++,
                  If[(a*d-b*c==em) && EvenQ[c] &&
                        OddQ[a] && (a*b*c*d<0),
                      fsm=fsm+(a*b/Abs[a*b])*(a*X+b)^n*(c*X+d)^tn;
                      ];
                  ];
              ];
          ];
      ];
    fsm=fsm/2;
    ssm=0;
    For[a=1,a<=em,a++,
      If[(Mod[em,a]==0) && OddQ[a],
          d=em/a;
          ssm=ssm+(a^n)*(nn^tn)*(X^w)*ber0[tn+1,d/(nn*X)]/(tn+1);
          ssm=ssm-(d^tn)*ber0[n+1,a*X]/(n+1);
          ];
      ];
    If[(j==2) && (Mod[m,nn]==0),
      For[d=1,d<=nn,d++,
        If[Mod[nn,d]==0,
          For[c=1,c<=m/nn,c++,
            If[Mod[m/nn,c]==0,
              ssm=ssm-(nn^w)*(X^w)*MoebiusMu[nn/d]*(c^tn)
                    *ber0[n+1,m*d/(c*nn^2*X)]/(d^n*(n+1));
              ];
            ];
          ];
        ];
      ];
    smwn[[j,i]]=fsm+ssm;
    cosmwn[[j,i]]=CoefficientList[Expand[smwn[[j,i]]],X];
    ];
  ];
m1=Array[b,{dw,dw}];
m2=Array[c,{dw,dw}];
For[i=1,i<=dw,i++,
    For[j=1,j<=dw,j++,m1[[i,j]]=cosmwn[[1,i]].cosmwn[[1,j]];
      m2[[i,j]]=cosmwn[[1,i]].cosmwn[[2,j]]]];
t=Inverse[m1].m2;
Print["representation matrix=",MatrixForm[t]];
Print["characteristic polynomial=",Det[t-x*IdentityMatrix[dw]]];

\end{verbatim}

The reader might be interested in comparing these results with those given at
``The Modular Form Database''
(http://modular.fas.harvard.edu/Tables/tables.html) by
William Stein.

\end{document}